\documentclass[twoside]{irmaems}
\usepackage{amssymb} 
\usepackage{amsmath} 
\usepackage{latexsym}
\usepackage[latin1]{inputenc}  
\usepackage[mathcal]{eucal}
\usepackage{psfrag}
\usepackage{graphicx}
\usepackage{makeidx}

\renewcommand\theequation{\arabic{section}.\arabic{equation}}

\theoremstyle{definition} 

\theoremstyle{plain}      

\makeindex

\begin{document}

\title{On Klein's \emph{So-called Non-Euclidean geometry}}
 

\author{Norbert A'Campo and Athanase Papadopoulos\thanks{The authors would like to thank Jeremy Gray and Fran\c cois Laudenbach for remarks on a preliminary version of this paper, and especially Hubert Goenner for his careful reading and for many interesting notes.
 The authors are supported by  the French ANR project FINSLER and the GEAR network of the National Science Foundation. Part of this work was done during a stay of the two authors at the Erwin Schr\"odinger Institute (Vienna).  The second author was also supported by a Tubitak 2221 grant  during a visit to Galatasaray University (Istanbul).}}

\address{Universit\"at Basel,  Mathematisches Institut, 
\\
Rheinsprung 21, 4051 Basel, Switzerland; 
\\
and 
Erwin Schr\"odinger International Institute of Mathematical Physics, 
\\
Boltzmanngasse 9, 1090, Wien, Austria
\\
email:\,\tt{Norbert.ACampo@unibas.ch}
\\
{}
\\
Institut de Recherche Math\'ematique Avanc\'ee,
\\
Universit{\'e} de Strasbourg and CNRS,\\
7 rue Ren\'e Descartes, 67084 Strasbourg Cedex, France;
\\
  Erwin Schr\"odinger International Inititute of Mathematical Physics, 
  \\
   Boltzmanngasse 9, 1090, Wien, Austria 
\\
and Galatasaray University, Mathematics Department,
\\
Ciragan Cad. No. 36, 
Besiktas, Istanbul, Turkey;
\\
email:\,\tt{papadop@math.unistra.fr}}

\maketitle
     
\begin{abstract}
 In two papers titled \emph{On the so-called non-Euclidean geometry, I} and \emph{II} (\cite{Klein-Ueber} and \cite{Klein-Ueber1}),  Felix Klein proposed a construction of the spaces of constant curvature -1, 0 and and 1  (that is, hyperbolic,  Euclidean and spherical geometry) within the realm of projective geometry. Klein's work was inspired by  ideas of Cayley who derived the distance between two points and the angle between two planes in terms of an arbitrary fixed conic in projective space. We comment on these two papers of Klein and we make relations with other works.
\end{abstract}

\bigskip 


\bigskip 


\bigskip

\tableofcontents

%

       \section{Introduction}
       
       Klein's model of hyperbolic space is well known to geometers. The underlying set is, in the planar case, the interior of an ellipse,  and in the three-dimensional case, the interior of an ellipsoid. The hyperbolic geodesics are represented by  Euclidean straight lines and  the distance between two distinct points  is defined as a constant times the logarithm of the cross ratio of the quadruple formed by this pair of points together with the two intersection points of the Euclidean straight line that joins them with the boundary of the ellipsoid, taken in the natural order.
       
      A much less known fact is that Klein also gave formulae for the distance in spherical and in Euclidean geometry using the cross ratio, taking instead of the ellipse (or ellipsoid) other kinds of conics. In the  case of Euclidean geometry, the conic is degenerate.\footnote{Recall that a conic is the intersection of a cone in 3-space and a plane. Degenerate cases occur, where the intersection is a single point (the vertex of the cone), or a straight ine (``counted twice"), or two intersecting straight lines.} In this way, the formulae that define the three geometries of constant curvature are of the same type, and the constructions of the three geometries are hereby done in a unified way in the realm of projective geometry. This is the main theme of Klein's two papers \emph{On the so-called non-Euclidean geometry, I} and \emph{II} (\cite{Klein-Ueber} \cite{Klein-Ueber1}) of Klein.
       
       In this paper, we review and comment on these two papers. 
       
       Klein's construction  was motivated by an idea of Cayley.\index{Cayley (Arthur)} We shall recall and explain this idea of Cayley and its developments by Klein. Let us note right away that although Klein borrowed this idea from Cayley, he followed, in using it, a different path. He even interpreted Calyey's idea in a manner so different from what the latter had in mind that Cayley misunderstood what Klein was aiming for, and  thought that the latter was mistaken.\footnote{\label{f:Klein} Cayley did not understand Klein's claim that the cross ratio is independent of the Euclidean underlying geometry, and therefore he disagreed with Klein's assertion that his construction of the three geometries of constant curvature was based only on projective notions; cf. \S \ref{s:Cayley}.} In fact, Cayley's interest was not in non-Euclidean geometry; it was rather to show the  supremacy of projective geometry over Euclidean geometry, by producing the Euclidean metric by purely projective methods. We shall comment on this fact in \S \ref{s:second-p} below.
       
The two papers by Klein were written in 1871 and 1872, just before Klein wrote his \emph{Erlangen program} manifesto \cite{Klein-Erlangen},\footnote{The first paper carries the date (``Handed in") ``D\"usseldorf, 19. 8. 1871" and the second paper: ``G\"ottingen, 8. 6. 1872". Thus the second paper was finished four months before the Erlangen Program, which carries the date ``October 1872".} the famous text in which he  he proposes a unification of all geometries based on the idea that geometry should be thought of as a transformation group rather than a space. Although this point of view is familiar to us, and seems natural today, this was not so for mathematicians even by the end of the nineteenth century.\footnote{In the introduction of his \emph{Lezioni di geometria proiettiva} \cite{Enriques-Lezioni}, Enriques writes that geometry studies the notion of space and the relation between its elements (points, curves, surfaces, lines, planes, etc.)} Without entering into the details of this philosophical question, let us recall that from the times of Euclid and until the raise of projective geometry, mathematicians were  reluctant to the use of the idea of transformation -- which, classically, carried the name  \index{motion}\emph{motion}\footnote{The word \emph{motion}, denoting a rigid transformation, was used by Peano. Hilbert used the word \emph{congruence}.}  -- as an element in the proof of a geometrical proposition.\footnote{The Arabic mathematician Ibn al-Haytham,\index{Ibn al-Haytham} (d. after 1040), in his book titled \emph{On the known} developed the first geometrical Euclidean system in which the notion of motion is a primitive notion (see \cite{Rashed} p. 446). Several centuries later, Pasch, Veronese and Hilbert  came up with the same idea. In his book \emph{La science et l'hypoth\`ese}, Poincar\'e discusses the importance of the notion of motion (see \cite{PK} p. 60).}

Klein's two papers 
 \cite{Klein-Ueber} and \cite{Klein-Ueber1} are actually referred to in the \emph{Erlangen program} text. We shall quote below some of Klein's statements from his \emph{Erlangen program} that are very similar to statements that are made in the two papers with which we are concerned. 

Another major element in the Erlangen program is the question of finding a \emph{classification} of the various existing geometries using the setting of projective geometry and of the projective transformation groups. Klein's paper \cite{Klein-Ueber} constitutes a leading writing on that subject, and it puts at the forefront of geometry both notions of transformation groups and of projective geometry. At the same time, this paper constitutes an important piece of work in the world of non-Euclidean geometry, and historically, it is probably the most important one in this domain, after the writings of the three founders of hyperbolic geometry (Lobachevsky, Bolyai and Gauss) and after  Beltrami's paper \cite{Beltrami-Saggio} on which we shall also comment below. 

In the introduction of the first paper \cite{Klein-Ueber}, Klein states that among all the works that were done in the preceding fifty years in the field of geometry, the development of projective geometry occupies the first place.\footnote{A similar statement is made in the introduction of the \emph{Erlangen program} \cite{Klein-Erlangen}: ``Among the advances of the last fifty years in the field of geometry, the development of \emph{projective} geometry occupies the first place".} Let us note that the use of the notion of \emph{transformation} and of \emph{ projective invariant} by geometers like Poncelet\footnote{Poncelet, for instance, made heavy use of  projective transformations in order to reduce proofs of general projective geometry statements to proofs of statements in special cases which are simpler to handle.} had prepared the ground for Klein's general idea that a geometry is a transformation group. The fact that the three constant curvature geometries (hyperbolic, Euclidean and spherical) can be developed in the realm of projective geometry is expressed by the fact that the transformation groups of these geometries are subgroups of the transformation group of projective transformations. Klein digs further this idea, namely, he gives explicit constructions of distances and of measures for angles\index{Klein's measure for distances}\index{Klein's measure for angles} in the three geometries\footnote{The reader will easily see  that it is a natural idea to define the notion of ``angle" at any point in the plane by using a conic (say a circle) at infinity, by taking the distance between two rays starting at a point as the length of the arc of the ellipse at infinity that the two rays contain However, this is not the definition used by Klein. His definition uses the cross ratio, like for distances between points, and this makes te result projectively invariant.} using the notion of ``projective measure" which was introduced by Cayley about twelve years before him.

 In fact, Cayley\index{Cayley (Arthur)} gave a construction of the Euclidean plane, equipped with its metric, as a subset of projective space, using projective notions.  This result is rather surprising because \emph{a priori} projective geometry is wider than Euclidean geometry insofar as the latter considers lines, projections and other notions of Euclidean geometry but without any notion of measurement of angles or of distances between points. Introducing distances between points or angle measurement and making the transformation group of Euclidean geometry a subgroup of the projective transformation group amounts to considering that space as a particular case of projective geometry; this was an idea of Cayley and,  before him, related ideas were emitted by Laguerre, Chasles and possibly others. However, Cayley did not use the notion of cross ratio in his definition of the distance. Klein's definitions of both measures (distances and angles) are based on the cross ratio. It was also Klein's contribution that the two non-Euclidean geometries are also special cases of a projective geometry.

      Klein starts his paper  \cite{Klein-Ueber} by referring to the work of Cayley, from which, he says, ``one may construct a projective measure on ordinary space using a second degree surface as the so-called fundamental surface". This sentence needs a little explanation. ``Ordinary space" is three-dimensional projective space. A ``measure" is a way of measuring distances between points as well as angles between lines (in dimension two) or between planes (in dimension three). Such a measure is said to be ``projective" if its definition is based on projective notions and if it is invariant under the projective transformations that preserve a so-called fundamental surface. Finally, the ``fundamental surface" is a quadric, that is, a second-degree surface, which is chosen as a ``surface at infinity" in projective space. 
      
      To say things briefly, a fixed quadric is  chosen. To define the distance between two points, consider the line that joins them; it intersects the quadric in two points (which may be real -- distinct or coincident -- or imaginary). The distance between the two points in the space is then, up to a constant factor,  the logarithm of the cross ratio of the quadruple formed by these points together with the intersection points of the line with the quadric taken in some natural order. The cross ratio (or its logarithm) could be imaginary, and the mutiplicative constant is chosen so that the result is real.   We shall come back to this definition in \S \ref{s:measures} of this paper. In any case, the projective measure depends on the choice of a fundamental surface and the definitions of measures on lines or on planes use dual constructions.  
Thus, Klein's construction is based on the fact that two points in the real projective space define a real line, which is also contained in a complex line (its complexification). If the real line  intersects the conic in two points, then these two points are real, and in this case the cross ratio is real. In the general case, the complex line intersects the conic in two points, which may be real or complex conjugate or coincident, and the cross ratio is a complex number. The multiplicative constant in front of the logarithm in the definition of the distance makes all distances real.

In some cases, there is a restriction on the conic; in the case where it is defined by the quadratic form $x^2+y^2-z^2$, the conic has an interior and an exterior,\footnote{A point is in the interior of the conic if there is not real tangent line from that point to the conic (a point that intersects it in exactly one point). Note that this notion applies only to real conics, since in the case of a complex conic, from a any point in the plane one can find a line which is tangent to the conic. This is expressed by the fact that a quadratic equation has a unique solution.} and one takes as underlying space the \emph{interior} of the conic.

Let us recall that in the projective plane, there are only two kinds of non-degenerate conics, viz. the \emph{real conics},\index{conic!real} which in homogeneous coordinates can be written as $x^2+y^2-z^2=0$, and the \emph{imaginary conics},\index{conic!imaginary}  which can be written as $x^2+y^2+z^2=0$. There is also a degenerate case where the conic is reduced to two coincident lines, which can be written in homogeneous coordinates as $z^2=0$, or also $x^2+y^2=0$. (Notice that this is degenerate because the differential of the implicit equation is zero.) In the way Cayley uses it, the degenerate conic can be thought of as the two points on the circle at infinity whose homogeneous coordinates are $(1,i,0)$ and  $(1,-i,0)$. 

Klein's work is based on Cayley's ideas of working with a conic at infinity which he called the ``absolute"\index{absolute conic} and measuring distances using this conic. Cayley gave a general formula that does not distinguish between the cases where the conic is real or imaginary, but he notes that Euclidean geometry is obtained in the case where the absolute degenerates into a pair of points. Klein makes a clear distinction between the cases of a real and an imaginary conic and he obtains the three geometries:
\begin{itemize}

\item The elliptic, in the case where the absolute is imaginary. (Notice that in this case, all the real directions are points in the projective space, since none of them intersects the imaginary conic.)  The fundamental conic in this case can be taken to be the imaginary circle whose equation is $\sum_i x_i^2=0$ (it has no real solutions).

\item The hyperbolic, in the case where the absolute is real. In this case, the conic has a well-defined ``interior" and an ``exterior", and the hyperbolic plane corresponds to the interior of the conic.

\item The parabolic, in the case where the absolute degenerates into two imaginary points. This is a limiting case of the preceding ones, and it corresponds to Euclidean geometry.
\end{itemize}
In this way, the three geometries become a particular case of projective geometry in the sense that the transformation groups of each geometry is a subgroup of the projective transformation group, namely the group of transformations that fix the given conic.

\section{Projective geometry}\index{geometry!projective}

As we already noted several times, in Klein's program, projective geometry acts as a unifying setting for many geometries. In fact, several theorems in Euclidean geometry (the theorems of Pappus, Pascal, Desargues, Menelaus,  Ceva, etc. ) find their real explanation in the setting of projective geometry. Let us say a few words of introduction to this geometry, since it will the main setting for what follows.

 For a beginner, projective geometry is, compared to the Euclidean,  a mysterious geometry. There are several reasons for that. First of all, the non-necessity of any notion of distance or of length may be misleading (What do we measure in this geo-metry?) Secondly, the Euclidean coordinates are replaced by the less intuitive (although more symmetric) ``homogeneous coordinates". Thirdly, in this geometry, lines intersect ``at infinity". There are points at infinity, there are ``imaginary points", and there is an overwhelming presence of the cross ratio, which, although a beautiful object, is not easy to handle. We can also add the fact that the projective plane is non-orientable and is therefore more difficult to visualize than the Euclidean.  One more difficulty is due to the fact that several among the founders of the subject had their particular point of view, and they had different opinions of what the fundamental notions should be. 
 
 Using modern notation, the ambient space for this geometry is the $n$-dimensional projective space $\mathbb{RP}^n$, that is, the quotient of Euclidean space $\mathbb{R}^{n+1}\setminus \{0\}$ by the equivalence relation which identifies a point $x$ with any other point $\lambda x$ for $\lambda \in \mathbb{R}^*$. The projective transformations of $\mathbb{RP}^n$ are  quotients of linear transformations of $\mathbb{R}^{n+1}$. They map lines, planes, etc.  in $\mathbb{R}^{n+1}$ to lines, planes, etc.  in $\mathbb{R}^{n+1}$; therefore, they map points, lines, etc.  in $\mathbb{RP}^n$ to points, lines, etc.  in $\mathbb{RP}^n$.  The incidence properties -- intersections of lines, of planes, alignment of points,  etc. -- are preserved by the projective transformations.
These transformations form a group called the projective linear group, denoted by $\mathrm{PGL}(n,\mathbb{R})$.  There is no metric on $\mathbb{RP}^n$ which is invariant by the action of this group, since this action is transitive on pairs of distinct points. As we already noted, Cayley\index{Cayley (Arthur)} observed that if we fix an appropriate quadric in $\mathbb{RP}^n$,  which he called the \emph{absolute}, we can recover the group of Euclidean geometry by restricting $\mathrm{PGL}(n,\mathbb{R})$ to the group of projective transformations that preserve this quadric. In this way, Euclidean space sits as the complement of the quadric, which becomes the \emph{quadric at infinity}.\index{quadric at infinity}
Klein states in his \emph{Erlangen program} \cite{Klein-Erlangen}: 
\begin{quote}\small Although it seemed at first sight as if the so-called metrical relations were not accessible to this treatment, as they do not remain unchanged by projection, we have nevertheless learned recently to regard them also from the projective point of view, so that the projective method now embraces the whole of geometry.
\end{quote}
Cayley's paper \cite{Cayley1859} on this subject was published in 1859 and it is abundantly cited by Klein in the two papers which are our main object of interest here. 

There is a situation which is familiar to any student in geometry, which is in the same spirit as Cayley's remark. If we fix a hyperplane $H$ in the projective plane $\mathbb{RP}^n$, then the subgroup of the group of projective transformations of $\mathbb{RP}^n$ that preserve $H$ (that is, that take this hyperplane to itself) is the affine group. Affine space is the complement of that hyperplane acted upon by the affine group. It is in this sense that ``affine geometry is part of projective geometry". In the projective space, at infinity of the affine plane stands a hyperplane. From Klein's point of view, affine geometry is determined by (and in fact it is identified with) the group of affine transformations, and this group is a subgroup of the group of projective transformations. Likewise, Euclidean and hyperbolic geometries are all part of affine geometry (and, by extension, of projective geometry), and furthermore, we have models of the spherical, Euclidean and hyperbolic spaces that sit in affine space, each of them with its metric and with a ``conic at infinity". One consequence is that all the theorems of projective geometry hold in these three classical geometries of constant curvature, and Klein insists on this fact, when he declares that projective geometry is ``independent of the parallel postulate" (see \S \ref{s:second-p} below).

In projective geometry, one studies properties of figures and of maps arising from projections (``shadows") and sections, or, rather, properties that are preserved by such maps. For instance, in projective geometry, a circle is equivalent to an ellipse (or to any other conic), since these objects can be obtained from each other by projection.

The first mathematical results where projective geometry notions are used, including duality theory, are contained in the works of Menelaus\index{Menelaus}, Ptolemy\index{Ptolemy} and in the later works of their Arabic commentators; see \cite{R1} and \cite{R2}. The Renaissance artists also used heavily projective geometry. A good instance of how projective geometry may be useful in perspective drawing is provided by Desargues' Theorem,\index{Theorem!Desargues} which is one of the central theorems of projective geometry and which we recall now. 

Consider in the projective plane two triangles $abc$ and $ABC$. We say that they are \emph{in axial perspectivity} if the three intersection points of lines $ab\cap AB, ac\cap AC, bc\cap BC$ are on a common line. We say that the three triangles are \emph{in central perspectivity} if the three lines $Aa, Bb, Cc$ meet in a common point. Desargues' theorem says that for any two triangles, being in axial perspectivity is equivalent to being in central perspectivity. The theorem is quoted in several treatises of perspective drawing. \footnote{Desargues' theorem was published for the first time by A. Bosse, in his \emph{Mani\`ere universelle de M. Desargues pour manier la perspective par petit pied comme le g\'eom\'etral} (Paris, 1648, p. 340). Bosse's memoir is reproduced in Desargues' \emph{\OE uvres} (ed. N. G. Poudra, Paris 1884, p. 413--415). Desargues' proof of the theorem uses Menelaus' Theorem. Von Staudt, in his \emph{Geometrie der Lage}, (Nuremberg, 1847) gave a proof of this theorem that uses only projective geometry notions.}

Ideas and constructions of projective geometry were extensively used by Renaissance artists like Leon Battista Alberti (1404-1472),\index{Alberti (Leon Battista)} Leonardo da Vinci\index{Da Vinci (Leonardo)} (1452-1519) and Albrecht D\"urer\index{D\"urer (Albrecht)} whom we already mentioned. All these artists used for instance the principle saying that any set of parallel lines in the space represented by the drawing which are not parallel to the plane of the picture must converge to a common point, called the \emph{vanishing point}.\footnote{The English term \emph{vanishing point} was introduced by Brook in the treatise \cite{Taylor} that he wrote in 1719, which is also the first book written in English on the art of perspective. The italian expression \emph{punto di fuga} was already used by Alberti and the other italian Renaissance artists.}  The italian mathematician and astronomer Guidobaldo del Monte\index{del Monte (Guidobaldo)} (1545 -1607) wrote a treaitise in six books (\emph{Perspectivae libri VI}, published in Pisa in 1600) in which he led the mathematical foundations of perspective drawing that included the vanishing point principle. In this treatise, the author often refers to Euclid's \emph{Elements}. The architect and famous scenographer\footnote{Designer of theatrical scenery} Nicola Sabbatini\index{Sabbatini (Nicola)} (1574-1654) made extensive use of Guidobaldo's theoretical work.  Guidobaldo del Monte's book is regarded as a mathematical work on a topic in projective geometry.

 It is usually considered that the modern theory of projective geometry started with  J.-V. Poncelet (1788-1867),\index{Poncelet (Jean-Victor)} in particular with his two papers, \emph{Essai sur les propri\'et\'es projectives des sections coniques} (presented at the French Academy of Sciences in 1820) and \emph{Trait\'e des propri\'et\'es projectives des figures} (1822). Poncelet tried to eliminate the use of coordinates and to replace them by synthetic reasonings. He made heavy use of duality (also called polarity) theory. This is based on the simple observation that in the projective plane any two points define a line and any two lines define a point. Using this fact, certain statements in projective geometry remain true if we exchange the words ``line" and ``point". A well known  example is Menelaus' Theorem which transforms under duality into Ceva's Theorem. Duality in projective geometry is at the basis of other duality theories in mathematics, for instance, the one in linear algebra, between a finite-dimensional vector space and the vector space of linear forms. To Poncelet is also attributed the so-called \emph{principle of continuity}\index{principle of continuity} which roughly says that the projective properties of a figure are preserved when the figure attains a limiting position. This permits him for instance to assert that points or lines, which disappear at infinity, become imaginary and can therefore be recovered, and one can make appropriate statements about them. 
 
  Among the other founders of projective geometry, we mention J. Brianchon (1783-1864),  A. F. M\"obius (1790-1868), M. Chasles\index{Chasles (Michel)} (1793-1880), K. G. K. von Staudt\index{von Staudt (Karl Georg Christian)} (1798-1867) and J. Steiner (1796-1863). We shall refer to some of them in the text below. The G\"ottingen lecture notes of Klein \cite{Klein-class} (1889-90) contain notes on the history of projective geometry (p. 61 \& ff.). We also refer the reader to the survey \cite{Enriques} by Enriques, in which the works of Cayley and of Klein are also analyzed. A concise modern historical introduction to projective geometry is contained in Gray's \emph{Worlds out of nothing} \cite{Gray}.

Perhaps Poncelet's\index{Poncelet (Jean-Victor)} major contribution, besides the systematic use of polarity theory,  was to build a projective geometry which is free from the analytic setting of his immediate predecessors and of the cross ratio (which he called the \emph{anharmonic ratio}, and so does Klein in the papers under consideration).
Chasles,\index{Chasles (Michel)} in his 1837 essay \emph{Aper\c cu historique sur l'origine et le d\'eveloppement des m\'ethodes en g\'eom\'etrie}, highlights the role of transformations in geometry, in particular in projective geometry, making a clear distinction between the metric and the projective (which he called ``descriptive") properties of figures. Von Staudt\index{von Staudt (Karl Georg Christian)}  insisted on the axiomatic point of view, and he also tried to build projective geometry independently from the notions of length and angle. One should also mention the work of E. Laguerre\index{Laguerre (Edmond)} (1843-1886), who was a student of Chasles and who, before Cayley, tried to develop the notions of Euclidean angle and distance relatively to a conic in the plane, cf. \cite{Laguerre} p. 66. Laguerre gave a formula for angle measure that involves the cross ratio. It is important nevertheless to note that Laguerre, unlike Klein, did not consider this as a possible \emph{definition} of angle. Laguerre's formula originates in the following problem that he solves: Given a system of angles $A,B,C,\ldots$ of a certain figure $F$ in a plane, satisfying an equation 
\[F(A,B,C,\ldots),\]
 find a relation satisfied by the image angles $A',B,'C',\ldots$ when the figure $F$ is transformed by a projective transformation (which Laguerre calls a homography). The solution that Laguerre gives is that $A',B,'C',\ldots$ satisfy the relation 
 \[F\left(
 \frac{\log a}{2\sqrt{-1}}, \frac{\log b}{2\sqrt{-1}},\frac{\log c}{2\sqrt{-1}},\ldots
 \right),
 \]
where $a,b,c,\ldots$ is the cross ratio of the quadruple of lines made up by two sides of the angles $A,B,C,\ldots$ together the lines $AP, AQ, BP, BQ, CP, CQp, \ldots$ which are the images of the lines made up by $A,B,C,\ldots$ and the two cyclic points of the plane of $F$.
Laguerre notes that Chasles, in his \emph{Trait\'e de g\'eom\'etrie sup\'erieure} \cite{Chasles-Traite}, p. 446, gave a solution to this problem in the case where the angles $A,B,C,\ldots$ share the same vertex or when they are equal. It seems that neither Cayley nor Klein, at the beginning of their work on this subject, were aware of Laguerre's work. Klein mentioned later on Laguerre's work, namely, in his 1889-90 lecture notes (\cite{Klein-class}, p. 47 and 61). In the \emph{Gesammelte  Mathematische Abhandlungen} (\cite{Klein-Ges}, vol. 1, p. 242) Klein declares that at the time he wrote his paper \cite{Klein-Ueber}, he was not aware of Laguerre's ideas. This work of Laguerre is also mentioned in \cite{Darboux} and \cite{Rouche}. See also \cite{Kl1926} for notes of Klein on Laguerre's work, and \cite{Troyanov} for some comments on Laguerre's formulae.

For Klein, the subject of projective geometry includes both its synthetic and its analytic aspect.\footnote{We can quote here Klein, from his \emph{Erlangen program} \cite{Klein-Erlangen}: ``The distinction between modern synthetic and modern analytic geometry must no longer be regarded as essential, inasmuch as both subject-matter and methods of reasoning have gradually taken a similar form in both. We choose therefore in the text as a common designation of them both the term \emph{projective geometry}".} In his historical remarks contained in his lecture notes (\cite{Klein-class} p. 61), he makes a distinction between the French and the German school of projective geometry and he notes that in the beginning of the 1850s, the French school had a serious advance over the German one, the latter still distinguishing between the projective and the metric properties.\footnote{In his \emph{Erlangen program} \cite{Klein-Erlangen}, Klein writes: ``Metrical properties are to be considered as projective relations to a fundamental configuration, the circle at infinity" and he adds in a note: ``This view is to be regarded as one of the most brilliant achievements of [the French school]; for it is precisely what provides a sound foundation for that distinction between properties of position and metrical properties, which furnishes a most desirable starting-point for projective geometry." Regarding Klein's comments on the difference between the French and the German schools, one may remember the context of that time, namely the French-German war (July 19, 1870-January 29, 1871), opposing the French Second Empire to the Prussian Kingdom and its allies; France suffered a crushing defeat and lost the Alsace-Moselle, which became the German Reichsland Elsa\ss -Lothringen.} 

 It is also fair to recall that this nineteenth-century activity on projective geometry was preceded by works of the Greeks, in particular by the work of Apollonius on the Conics (where the notion of a polar line with respect to a conic appears for the first time),  by works of Pappus, and by the much later works of several French mathematicians, including G. Desargues (1591-1661), who wanted to give a firm mathematical basis to the perspective theory used by painters and architects, and then, B. Pascal (1623-1662), who was influenced by Desargues, and G. Monge (1746-1818). 
 
 \section{Non-Euclidean geometry}
 
In this section, we recall a few facts on the birth and the reception of non-Euclidean geometry, and on its relation with projective geometry.

 Nikolai Ivanovich Lobachevsky\index{Lobachevsky (Nikolai Ivanovich)} was the first to publish a treatise on hyperbolic geometry, namely, his \emph{Elements of Geometry} \cite{Loba-Elements} (1829).  The two other founders of the subject are J\'anos Bolyai\index{Bolyai (J\'anos)} and Carl Friedrich Gauss.\index{Gauss (Carl Friedrich)} For more than 50 years,  these three works remained unknown to the mathematical community. Lobachevsky's work was acknowledged as being sound only ten years after his death (1856), after Gauss's correspondence was published.\footnote{In a letter he wrote to his friend the astronomer H. C. Schumacher, dated 28 November 1846, Gauss expresses his praise for Lobachevsky's work, cf. \cite{Gauss} p. 231--240, and it was after the publication of this letter that mathematicians started reading Lobachevsky's works. Lobachevsky was never aware of that letter. Gauss's correspondence was published during the few years that followed Gauss's death, namely, between 1860 and 1865.}
This geometry first attracted the attention of Cayley\index{Cayley (Arthur)} and then Beltrami. Beltrami\index{Beltrami (Eugenio)} started by publishing two papers on the subject, \emph{Saggio di Interpretazione della geometria non-Euclidea} \cite{Beltrami-Saggio} (1868) which concerns the two-dimensional case and  \emph{Teoria fondamentale degli spazii di curvatura costante}  \cite{Beltrami-Teoria} (1868-69) which concerns the three-dimensional case.\footnote{Eugenio Beltrami\index{Beltrami (Eugenio)} (1835-1900) was born in a family of artists.  He spent his childhood in a period of political turbulence:  the Italian revolutions, the independence war, and eventually the unification of Italy. He studied mathematics in Pavia between 1853 and 1856, where he followed the courses of Francesco Brioschi, but due to lack of money or may be for  other reasons (Loria reports that Beltrami was expelled from the university because he was accused of promoting disorders against the rector \cite{Loria}), he interrupted his studies and took the job of secretary of the director of the railway company in Verona. The first mathematical paper of Beltrami was published in 1862, and in the same year he got a position at the University of Bologna. He later on moved between several universities, partly because of the political events in Italy, and he spent his last years at the university of Rome. A stay in Pisa, from 1863 to 1866, was probably decisive for his mathematical future research; he met there Betti and Riemann (who was in Italy for health reasons). Two of the most influential papers of Beltrami are quoted in the present  survey,  \cite{Beltrami-Saggio} (1868) and \cite{Beltrami-Teoria} (1869). They were written during his second stay in Bologna where he was appointed on the chair of rational mechanics.  His name is attached to the Beltrami equation, a fundamental equation in the theory of quasiconformal mappings, and to the Laplace-Beltrami operator.  Besides mathematics, Beltrami cultivated physics, in particular thermodynamics, fluid dynamics, electricity and magnetism. He  translated into Italian the work of Gauss on conformal representations. He contributed to the history of mathematics by publishing a paper on the work of Saccheri on the problem of parallels (\emph{Un precursore italiano di Legendre e di Lobatschewski}, 1889), comparing this work to the works of Borelli, Clavius, Wallis, Lobachevsky and Bolyai on the same subject, and highlighting the results on non-Euclidean geometry that are inherent in that work. Besides mathematics, Beltrami cultivated music, and also politics.  In 1899, he became (like his former teacher Brioschi) senator of the Kingdom of Italy.}
We shall elaborate on them below. 

Lobachevsky, Bolyai and Gauss developed the hyperbolic geometry system starting from the axiomatic point of view, that is, drawing conclusions from the set of axioms of Euclidean geometry with the parallel axiom replaced by its negation. Beltrami,\index{Beltrami (Eugenio)} in his paper \cite{Beltrami-Saggio},  was the first to establish the relation between hyperbolic geometry and negative curvature. None of the three founders of hyperbolic geometry used the notion of curvature. It is true that at that time, curvature in the Riemannian geometry sense was not yet discovered, but Gauss had already introduced the notion of surface curvature and he had showed that it is independent of any embedding in an ambient Euclidean space. Gauss nevertheless did not make explicit the relation between curvature and the geometry of the hyperbolic plane.  

The definition of spherical geometry as a system which can also be defined using the notion of constant positive curvature is due to Riemann. It is also a geometrical system which is at the same level as Euclidean geometry, where the ``lines"  are the great circles of the sphere but where there are no disjoint lines. It is also good to recall that, unlike hyperbolic geometry, the geometrical system of the sphere cannot be obtained from the Euclidean one by modifying only one axiom, since not only the Euclidean parallel axiom is not valid on the sphere, but other axioms as well, e.g. the one saying that lines can be extended indefinitely.\footnote{The intuition that there are exactly \emph{three} geometries, and that these three geometries are the hyperbolic, the Euclidean and the spherical, can be traced back to older works. In the memoir \emph{Theorie der Parallellinien}  of J. H. Lambert \cite{Lambert}, written in 1766, that is, more than 100 years before Klein wrote his memoir \cite{Klein-Ueber}, the author, attempting a proof of Euclid's parallel postulate, developed a detailed analysis of geometries that are based on three assumptions, concerning a class of quadrilaterals, which are now called \emph{Lambert} or  \emph{Ibn al-Haytham-Lambert quadrilaterals}. These are quadrilaterals having three right angles, and the assumptions he made are that the fourth angle is either acute, right or obtuse. These assumptions lead respectively to hyperbolic, Euclidean and spherical geometry. One must add that Lambert was not the first to make such a study of these quadrilaterals. Gerolamo Saccheri (1667-1733) and, before him,  Ab\=u `{}Al\=\i \  Ibn al-Haytham and `Umar al-Khayy\=am (1048-1131) made similar studies. Of course, in all these works, the existence of hyperbolic geometry was purely hypothetical. The approaches of these authors consisted in assuming that such a geometry exists and to try to deduce a contradiction. We refer the interested reader to the recent edition of Lambert's work \cite{PT}, with a French translation and mathematical comments.} Riemann established the bases of spherical geometry in his famous habilitation lecture  \emph{\"Uber die Hypothesen, welche der Geometrie zu Grunde liegen} (On the hypotheses which lie at the foundations of geometry)  (1854)  \cite{Riemann-hypothesis}. It is generally accepted that the three geometries -- Euclidean, hyperbolic and spherical -- appear clearly for the same time at the same level, as ``the" three geometries of constant curvature in the paper \cite{Klein-Ueber} of Klein. However, one can mention a letter from Ho\"uel to De Tilly, dated April 12, 1872, in which he writes:\footnote{A few words are needed on Ho\"uel and de Tilly, two major major figures in the history of non-Euclidean geometry but whose names remain rather unknown to most geometers. Guillaume-Jules Ho\"uel\index{Ho\"uel (Guillaume-Jules)} (1823--1886) taught at the University of Bordeaux, and he wrote geometric treatises giving a modern view on Euclid's {\it Elements}. He was working on the impossibility of proving the parallel postulate  when, in 1866, he came across the writings of Lobachevsky, and became convinced of their  correctness.  In the same year, he  translated into French Lobachevsky's   {\it Geometrische Untersuchungen zur Theorie der Parallellinien}  together with excerpts from the correspondence between Gauss and Schumacher on non-Euclidean geometry, and he published them in the {\it M\'emoires de la Soci\'et\'e des Sciences physiques et naturelles de Bordeaux}, a journal of which he was the editor. Ho\"uel also translated  into French and published in French and Italian  journals works by several other authors on non-Euclidean geometry, including Bolyai, Beltrami, Helmholtz, Riemann and Battaglini.   Barbarin, in his book 
   \emph{La G\'eom\'etrie non Euclidienne}  (\cite{Barbarin1928} p. 12)
 writes  that Ho\"uel, ``who had an amazing working force, did not hesitate to learn all the European languages in order to make available to his contemporaries the most remarkable mathematical works." Ho\"uel also solicited for his journal several papers on hyperbolic geometry, after the French Academy of Sciences, in the 1870s, decided to refuse to consider papers on that subject. We refer the interested reader to the article by Barbarin \cite{Barbarin-Correspondance} and the forthcoming edition of the correspondence between Ho\"uel and de Tilly \cite{Henry-Nabonnand-H}.
Beltrami had a great respect for Ho\"uel, and  
there is a very interesting correspondence between the two men, see  \cite{Beltrami-Boi}. It appears from these letters that Beltrami's famous \emph{Saggio di Interpretazione della geometria non-Euclidea} \cite{Beltrami-Saggio} arose from ideas that he got after reading Lobachevsky's {\it Geometrische Untersuchungen zur Theorie der Parallellinien}  in the French translation by Ho\"uel, see \cite{Beltrami-Boi} p. 9. For a detailed survey on the influence of Ho\"uel's work see \cite{Brunel}. 
 
 Joseph-Marie de Tilly\index{De Tilly (Joseph-Marie)} (1837--1906) was a member of the Royal Belgian Academy of Sciences, and he  was also an officer in the Belgian army, teaching mathematics at the Military School. In the 1860s, de Tilly, who was not aware of the work of Lobachevsky, developed independently a geometry in which Euclid's parallel postulate does not hold. One of his achievements is the introduction of the notion of distance as a primary notion in the three geometries: hyperbolic, Euclidean and spherical. He developed an axiomatic approach to these geometries based on metric notions, and he highlighted some particular metric relations between finite sets of points; see for instance his  {\it Essai sur les Principes Fondamentaux de la G\'eom\'etrie et de la M\'ecanique} \cite{Tilly-Principes} and his
 {\it Essai de G\'eom\'etrie analytique g\'en\'erale} \cite{Tilly-Analytique}.
} (see \cite{Henry-Nabonnand-H}): 
\begin{quote}\small
The idea of the three geometries is not due to Klein: it goes back to Lejeune-Dirichlet, who has thoroughly meditated upon this subject, but who, unfortunately, did not leave us anything written.
\end{quote}

In Klein's paper \cite{Klein-Ueber}, while the three geometries are placed at the same level of importance, Euclidean geometry acts as a \emph{transitional} geometry\index{geometry!transitional} between the other two. Klein writes about this:
\begin{quote}\small 
 Straight lines have no points at infinity, and indeed one cannot draw any parallel at all to a given line through a point outside it.
\\
A geometry based on these ideas could be placed alongside ordinary Euclidean geometry like the above-mentioned geometry of Gauss, Lobachevsky and Bolyai. While the latter gives each line two points at infinity, the former gives none at all (i.e. it gives two imaginary points). Between the two, Euclidean geometry stands as a transitional case; it gives each line two coincident points at infinity.
\end{quote}
 We develop this idea of ``transitional geometry" in our paper \cite{ACPK2} in this volume.

 Today, people are so much used to these ideas that it is hard for them to appreciate their novelty for that epoch and their importance. Let us recall in this respect that Klein's paper came out only three years after Beltrami published his two famous papers in which he confirmed that Lobachevsky's researches on hyperbolic geometry were sound.  
As we shall see later in this paper, to prove that Cayley's constructions lead to the non-Euclidean geometries, Klein essentially argued in an synthetic way, at the level of the axioms, showing that the characteristics of the Lobachevsky and of the spherical geometries are satisfied in the geometry defined by this distance function. But he also described the differential-geometric aspects, introducing a  notion of curvature which he showed is equivalent to Gauss's surface curvature.  
  
In the rest of this paper, we shall present the basic ideas contained in Klein's two papers, making connections with other ideas and works on the same subject.

\section{Preliminary remarks on Klein's papers}

In this section, we start by summarizing the important ideas contained in Klein's two  papers, and then we discuss the reception of these ideas by Klein's contemporaries and by other mathematicians.  We then make some remarks on the names \emph{hyperbolic, elliptic} and \emph{parabolic} geometries that were used by Klein.

\medskip

Klein's major contributions in these two papers include the following:
 
\begin{enumerate}
\item An explanation of the notion of Cayley measure\index{Cayley measure} and its representation, including it in two important settings: transformation groups and curvature.
\item A realization of Lobachevsky's geometry as a metric space (and not only as a system of axioms) where the metric is given by the Cayley measure.
\item The construction of a new model of Lobachevsky's geometry, by taking, in Cayley's construction, the ``absolute" to be an arbitrary real second-degree curve in the projective plane and showing that the interior of that curve, equipped with some adequate structure, is a model of Lobachevsky's geometry. 
Although the idea for the construction originates in Cayley's work (Cayley gave a formula for a distance function without realizing that the resulting metric space is the Lobachevsky space), and although the construction of such a model for the hyperbolic plane (but without the distance function) had been made three years earlier by Beltrami in this paper \cite{Beltrami-Saggio} (1868) in which he realized that the Euclidean segments of the disk are models for the geodesics of hyperbolic space, Klein gave the first explicit distance function for hyperbolic geometry. At the same time, this made the first  explicit link between hyperbolic geometry and projective geometry. 
\item A unified setting for Euclidean, hyperbolic and spherical geometries, as these three geometries can be considered as special cases of projective geometry. Although it is well known that Klein gave a formula for the hyperbolic metric using cross ratio, it is rather unknown to modern geometers that Klein also gave in the same way formulae for the elliptic and for the Euclidean distance functions using the cross ratio. Cayley expressed the advantage of Klein's distance formula in his comments on his paper \cite{Cayley1859} contained in his \emph{Collected mathematical papers} edition (\cite{Cayley-collected} Vol. II, p. 604):
\begin{quote}
In his first paper, Klein substitutes, for my $\cos^{-1}$ expression for the distance between two points,\footnote{The formulae to which Cayley refers are contained in his paper  \cite{Cayley1859} p. 584-585. After giving these formulae in the case where the absolute is a general conic, he writes:
\begin{quote}\small
The general formulae suffer no \emph{essential} modifications, but they are greatly simplified by taking for the point-equation of the absolute
\[x^2+y^2+z^2=0,\]
or, what is the same, for the line-equation
\[\xi^2+\eta^2+\zeta^2=0.\]
In fact, we then have for the expression of the distance of the points $(x,y,z), (x',y',z')$, 
\[
\cos^{-1}\frac{xx'+yy'+zz'}{\sqrt{x^2+y^2+z^2}\sqrt{x'^2+y'^2+z'^2}};\]
for that of the lines $(\xi,\eta,\zeta), (\xi',\eta',\zeta')$, 
\[
\cos^{-1}\frac{\xi\xi + \eta\eta'+\zeta\zeta'}{\sqrt{\xi^2+\eta^2+\zeta^2}\sqrt{\xi'^2+\eta'^2+\zeta'^2}};\]
and that for the point $(x,y,z)$ and the line $ (\xi',\eta',\zeta')$,
\[
\cos^{-1}\frac{\xi' x+\eta' y+\zeta' z}{\sqrt{x^2+y^2+z^2}\sqrt{\xi'^2+\eta'^2+\zeta'^2}}.\]
 \end{quote}
The reader will notice the analogy between these formulae and the familiar formula for distance in spherical geometry (the ``angular distance"), which, Klein also establishes in his paper \cite{Klein-Ueber}; see Formulae (\ref{FK1}) and  (\ref{FK2}) in the present paper.} a logarithmic one; viz. in linear geometry if the two fixed points are $A, B$ then the assumed definition for the distance of any two points $P,Q$ is
\[\mathrm{dist}. \, (PQ)=s\log \frac{AP.BQ}{AQ.BP};\]
this is a great improvement, for we at once see that the fundamental relation, $\mathrm{dist}. \, (PQ)+\mathrm{dist}. \, (QR) = \mathrm{dist}. \, (PR)$, is satisfied.
\end{quote}
We note that Cayley, instead of using the cross ratio in his definition of distances, showed that his formulae are invariant by the action of the projective geometry transformations on homogeneous coordinates.
\item Likewise, Klein gave a formula for the dihedral angle between two planes as a cross ratio between four planes, the additional two planes being the tangent planes to the fixed conic passing through the intersection of the first two planes. 
\item The conclusion that each of the three geometries is consistent if projective geometry is consistent.
\item The idea of a \emph{transitional geometry}, that is, a geometrical system in which one can transit continuously from spherical to hyperbolic geometry, passing through Euclidean geometry. 
  \item The introduction of the names hyperbolic, parabolic and elliptic for the Lobachevsky, Euclidean and spherical geometries respectively, thus making the relation with other settings where the three words ``hyperbolic", ``parabolic" and ``elliptic" were already used. We shall discuss this in the next section.
\end{enumerate}
We shall elaborate on all these items below.

Klein's papers are sometimes difficult to read and they were received by the mathematical community in diverse manners. Let us quote, for example, Darboux, from his obituary concerning Henri Poincar\'e \cite{Darboux}: 
\begin{quote}\small
Mr. Felix Klein is the one who removed these very serious objections [concerning non-Euclidean geometry] by showing in a beautiful memoir that a geometry invented by the famous Cayley and in which a conic called the absolute provides the elements of all measures and enables, in particular, to define the distance between two points, gives the most perfect and adequate representation of non-Euclidean geometry.
\end{quote}
On the other hand, Genocchi\footnote{Angelo Genocchi (1817-1889) was an Italian mathematician who made major contributions in number theory, integration and the theory of elliptic functions. Like Cayley, he worked for several years as a lawyer, and he taught law at the University of Piacenza, but at the same time he continued cultivating mathematics with passion. In 1859, he was appointed professor of mathematics at the University of Torino, and he remained there until 1886. During the academic year 1881-82, Guiseppe Peano served as his assistant, and he subsequently helped him with his teaching, when Genocchi became disabled after an accident. Genocchi's treatise \emph{Calcolo differenziale e principii di calcolo integrale con aggiunte del Dr. Giuseppe Peano}, written in 1864, was famous in the Italian universities.} wrote, regarding the same matter  (\cite{Genocchi} p. 385): 
\begin{quote}\small
From the geometric point of view, the spirit may be shocked by certain definitions adopted by Mr. Klein: the notions of distance and angle, which are so simple, are replaced by complicated definitions [...] The statements are extravagant. 
\end{quote}
Hans Freudenthal (1905-1990), talking about Klein's analysis of the work of von Staudt on the so-called ``fundamental theorem of projective geometry" in which Klein discusses some continuity issues that were missing in von Staudt's arguments,\footnote{This continuity issue is mentioned in Chapter 2 of this volume \cite{Gray-K}, and its is discussed in detail in \cite{Voelke2008}.} writes \cite{Freudenthal}:
\begin{quote}\small
[In 1873], logical analysis was not the strong point of Klein, and what he wrote on that question in the years that followed is as much confusing as possible.
\end{quote}

We end this section with two remarks. The first one concerns the title of the two papers  \cite{Klein-Ueber} and  \cite{Klein-Ueber1}, and the second concerns the names ``hyperbolic", elliptic" and ``parabolic" geometries.

Klein's title, \emph{On the so-called Non-Euclidean geometry}, may be considered as having a negative connotation, and indeed it does. This is also the title of a note (Note No. 5) at the end of his \emph{Erlangen program} text. In that note, Klein writes:
\begin{quote}\small We associate to the name Non-Euclidean geometry a crowd of ideas that have nothing mathematical, which are accepted on the one hand with as much enthusiasm that they provoke aversion on the other hand, ideas in which, in any case, our exclusively mathematical notions have nothing to do.
\end{quote}
However, Klein, in his later writings, used extensively the term ``non-Euclidean geometry", without the adjective ``so-called".

Now about the names of the three geometries.

Klein coined the expressions ``elliptic", ``hyperbolic" and ``parabolic" geometry\index{geometry!elliptic}\index{geometry!parabolic}\index{geometry!hyperbolic}  as alternative names for  spherical, Lobachevsky and Euclidean geometry respectively. In his \emph{\"Uber die sogenannte Nicht-Euklidische Geometrie} \cite{Klein-Ueber}, he writes (Stillwell's translation p. 72): 
\begin{quote}\small Following the usual terminology in modern geometry, these three geometries will be called \emph{hyperbolic, elliptic or parabolic} in what follows, according as the points at infinity of a line are real, imaginary or coincident.
\\
These three geometries will turn out to be special cases of the general Cayley measure.\index{Cayley measure} One obtains the parabolic (ordinary) geometry by letting the fundamental surface for the Cayley measure degenerate to an imaginary conic section. If one takes the fundamental surface to be a proper, but imaginary, surface of second degree, one obtains the elliptic geometry. Finally, the hyperbolic geometry is obtained when one takes the fundamental surface to be a real, but not ruled, surface of second degree and considers the points inside it.
\end{quote}

 In the paragraph that precedes this one, Klein makes the statement that in hyperbolic geometry, straight lines have two points at infinity, that in spherical geometry, they have no point at infinity, and that in Euclidean geometry, they have two coincident points at infinity. This can explain the meaning of the expression ``following the usual terminology", if we recall furthermore that in the projective plane, a hyperbola meets the line at infinity in two points, a parabola meets it is one point and an ellipse does not meet it at all. In a note in \cite{Sources}  (Note 33 on the same page), Stillwell recalls that before Klein, the points on a differentiable surface were called hyperbolic, elliptic or parabolic according to whether the principal tangents at these points are real, imaginary or coincident. He also recalls that Steiner used these names for certain surface involutions,  an involution being called hyperbolic, elliptic  or parabolic depending on whether the double points arising under it are respectively real, imaginary  or coincident. The name \emph{non-Euclidean geometry} is due to Gauss.\footnote{Gauss used it in his correspondence with Schumacher.}
 
Soon after Klein introduced this terminology, Paul du Bois-Reymond (1831-1889) introduced (in 1889) the classification of second-order differential operators into ``elliptic", ``hyperbolic" and ``parabolic". 

Let us note finally that it is easy to be confused concerning the order and the content of the two papers of Klein, if one looks at the French versions. The papers appeared in 1871 and 1873 respectively, under the titles \emph{\"Uber die sogenannte Nicht-Euklidische Geometrie} and 
 \emph{\"Uber die sogenannte Nicht-Euklidische Geometrie (Zweiter Aufsatz)}. In 1871, and before the first paper \cite{Klein-Ueber} was published, a short version, presented by Clebsch,\footnote{Alfred Clebsch (1833-1872) was a young professor at G\"ottingen,  who was responsible for Klein's first invitation at that university, in 1871. He was well aware of Cayley's work on invariant theory, and he transmitted it to Klein. Klein stayed in G\"ottingen a few months, and then moved to Erlangen, where he was appointed professor, again upon the recommendation of Clebsch. He came back to the University of G\"ottingen in 1886, and he stayed there until his retirement in 1913. Clebsch was also the founder of the \emph{Mathematische Annalen}, of which Klein became later on one of the main editors. See also \cite{Gray-K}.}, appeared in the \emph{Nachrichten von der Kgl. Gesellschaft der Wissenschaften zu G\"ottingen}, under the title \emph{F. Klein,  \"Uber die sogenannte Nicht-Euklidische Geometrie. Vorgelegt von A. Clebsch.} The same year, a translation of this short paper appeared under the title \emph{Sur la g\'eom\'etrie dite non euclidienne, de F\'elix Klein}, in the \emph{Bulletin de sciences math\'ematiques et astronomiques}, translated by Ho\"uel. A translation by Laugel of the first paper (1871) appeared much later in the \emph{M\'emoires de la Facult\'e des Sciences de Toulouse} under the title \emph{Sur la G\'eom\'etrie dite non euclidienne, par Mr. F\'elix Klein}, in 1898. In the volume \cite{Sources} (1996) which contains translations by J. Stillwell of some of the most important sources on non-Euclidean geometry,   only the first paper \cite{Klein-Ueber} by Klein is included, under the title \emph{On the so-called non-Euclidean geometry}, and it is followed by a short excerpt (6 lines) of the second paper.

             \section{The work of Cayley} \label{s:Cayley}
      In this section, we comment on the idea of Cayley\footnote{Arthur Cayley\index{Cayley (Arthur)} (1821-1895) was born in the family of an English merchant who was settled in Saint-Petersburg. The family returned to England when the young Arthur was eight. Cayley is one of the main inventors of the theory of invariants. These include invariants of algebraic forms (the determinant being an example), and algebraic invariants of geometric structures and the relations they satisfy (``syzygies"). Cayley studied mathematics and law. He was very talented as a student in mathematics, and he wrote several papers during his undergraduate studies, three of which were published in the \emph{Cambridge Mathematical Journal}. The subject included determinants, which became later on one of his favorite topics.  After completing a four-year position at Cambridge university, during which he wrote 28 papers for the Cambridge journal, Cayley did not succeed in getting a job in academics. He worked as a lawyer during 14 years, but he remained active mathematics; he wrote during these years about 250 mathematical papers. In 1863, he was appointed professor of mathematics at Cambridge. His list of papers includes about 900 entries, on all fields of mathematics of his epoch. The first definition of an abstract group is attributed to him, cf. his paper \cite{Cayley-group}. Cayley proved in particular that every finite group $G$ is isomorphic to a subgroup of a symmetric group on $G$. His name is attached to the famous \emph{Cayley graph} of a finitely generated group, an object which is at the basis of modern geometric group theory. Cayley is also one of the first discoverers of geometry in dimensions greater than 3. In his review of Cayley's \emph{Collected Mathematical Papers} edition in 13 volumes, G. B. Halsted writes: ```Cayley not only made additions to every important subject of pure mathematics, but whole new subjects, now of the most importance, owe their existence to him. It is said that he is actually now the author most frequently quoted in the living world of mathematicians" \cite{Halsted-Cayley}. We refer the reader to the biography by Crilly \cite{Crilly} which is regrettably short of mathematical detail, but otherwise very informative and accurate.} which acted as a motivation for Klein's work.

Let us quote again Klein, from the introduction to his paper \cite{Klein-Ueber}:
\begin{quote}\small
\emph{It is our purpose to present the mathematical results of these works [of Gauss, Lobachevsky and Bolyai], insofar as they relate to the theory of parallels, in a new and intuitive way, and to provide a clear general understanding.}
\\
The route to this goal is through projective geometry. By the results of Cayley, one may construct a projective measure on ordinary space using an arbitrary second degree surface as the so-called fundamental surface. Depending on the type of the second degree surface used, this measure will be a model for the various theories of parallels in the above-mentioned works. But it is not just a model for them; as we shall see, it precisely captures their inner nature.
\end{quote}

The paper to which Klein refers is Cayley's \emph{Sixth Memoir upon Quantics}\footnote{In Cayley's terminology, a quantic is a homogeneous polynomial.} \cite{Cayley1859} which appeared in 1859. In this paper, Cayley asserts that descriptive geometry (which is the name he used for projective geometry) ``is all geometry", an idea which was taken up by Klein later on.\footnote{In fact, the statement is also true if we interpret it in the following sense (which, however, is not what Cayley meant): Most of the work that was being done by geometers at the time Cayley made that statement was on projective geometry.} In particular, Cayley considered that projective geometry includes metrical geometry (which is the name he  used for Euclidean geometry) as a special case. In Cayley's words: ``A chief object of the present memoir is the establishment, upon purely descriptive principles, of the notion of distance". At first sight, there is something paradoxical in this statement, because length is not a projective notion. In fact, in his foundational work on descriptive geometry, and in particular in his  famous 1822 \emph{Trait\'e} \cite{Poncelet-Traite},  Poncelet had already stressed on the distinction between the metrical properties (namely, those that involve distance and angle), which are not preserved by projective transformations, and the projective (which he calls ``descriptive") properties, which are precisely the properties preserved by projective transformations, e.g. alignment of points, intersections of lines, etc. Thus, in principle, there are no distances, no circles and no angles in projective geometry. Cayley, followed by Klein, was able to define such notions using the concepts of projective geometry by fixing a quadric in projective space, in such a way that these properties are invariant under the projective transformations that fix the quadric. The cross ratio of four points is a projective invariant, and in some sense it is a complete projective invariant, since a transformation of projective space which preserves the cross ratio of quadruples of aligned points is a projective transformation. Therefore, it is natural to try to define distances and angles using the cross ratio. This is what Klein did. Likewise, it was an intriguing question, addressed by Klein, to try to express the concept of parallelism in Euclidean and in hyperbolic geometry using projective notions, although parallelism is a priori not part of projective geometry. Cayley defined a geometry which is non-Euclidean, but did not realize that it coincides with the Lobachevsky geometry. Let us quote Cayley's paper (the conclusion):
\begin{quote}\small
I have, in all that has preceded, given the analytical theory of distance along with the geometrical theory, as well for the purpose of illustration, as because it is important to have an analytical expression of a distance in terms of the coordinates; but I consider the geometrical theory as perfectly complete in itself: the general result is as follows; viz. assuming in the plane (or space of geometry of two dimensions) a conic termed the \emph{absolute}, we may by means of this conic, by descriptive constructions, divide any line or range of points whatever, and any point or pencil of lines whatever, into an infinite series of infinitesimal elements, which are (as a definition of distance) assumed to be equal; the number of elements between any two points of the range or two lines of the pencil, measures the distance between the two points or lines; and by means of the pencil, measures the distance between the two points or lines; and by means of the quadrant, as a distance which exists as well with respect to lines as points, we are enabled to compare the distance of two lines with that of two points; and the distance of a point and a line may be represented indifferently as the distance of two points, or as the distance of two lines.
\\
In ordinary spherical geometry, the theory undergoes no modification whatever; the absolute is an actual conic, the intersection of the sphere with the concentric evanescent sphere.
\\
In ordinary plane geometry, the absolute degenerates into a pair of points, viz. the points of intersection of the line at infinity with any evanescent circle, or what is the same thing, the absolute is the two circular points at infinity. The general theory is consequently modified, viz. there is not, as regards points, a distance such as the quadrant, and the distance of two lines cannot be in any way compared with the distance of two points; the distance of a point from a line can be only represented as a distance of two points.
\\
I remark in conclusion that, \emph{in my point of view}, the more systematic course in the present introductory memoir on the geometrical part of the subject of quantics, would have been to ignore altogether the notions of distance and metrical geometry; for the theory in effect is, that the metrical properties of a figure are not the properties of the figure considered \emph{per se} apart from everything else, but its properties when considered in connexion with another figure, viz. the conic termed the absolute. The original figure might comprise a conic; for instance, we might consider the properties of the figure formed by two or more conics, and we are then in the region of pure descriptive geometry by fixing upon a conic of the figure as a standard of reference and calling it the absolute. Metrical geometry is thus a part of descriptive geometry, and descriptive is \emph{all} geometry and reciprocally; and if this can be admitted, there is no ground for the consideration in an introductory memoir, of the special subject of metrical geometry; but as the notions of distance and of metrical geometry could not, without explanation, be thus ignored, it was necessary to refer to them in order to show that they are thus included in descriptive geometry.
\end{quote}

In his \emph{Lectures on the development of mathematics in the XIXth century} \cite{Klein19} (1926-1927),  Klein recounts how he came across Cayley's ideas (p. 151):
\begin{quote}\small 
In 1869, I had read Cayley's theory in the version of Fiedler\footnote{[Wilhelm Fiedler (1832-1911)]} of Salmon's \emph{Conics}. Then, I heard for the first time the names of Bolyai and Lobatscheffski, from Stolz,\footnote{Otto Stolz\index{Stolz (Otto)} (1842-1905) was a young mathematician at the time when Klein met him. He obtained his habilitation in Vienna in  1867 and, starting from 1869, he  studied in Berlin under Weierstrass, Kummer and Kronecker. He attended Klein's lecture in 1871 and he remained in contact with him. He became later on a successful textbook writer.} in the winter of 1869/70, in Berlin. From these indications I had understood very little things, but I immediately got the idea that both things should be related. In February 1870, I gave a talk at Weierstrass's seminar on the Cayley metric.\footnote{In 1870, Weierstrass started at the university of Berlin a seminar on non-Euclidean geometry.} In my conclusion, I asked whether there was a correspondence with Lobatscheffski. The answer I got was that these were two very different ways of thinking, and that for what concerns the foundations of geometry, one should start by considering the straight line as the shortest distance between two points. I was daunted by this negative attitude and this made me put aside the insight which I had. [...]
\\
In the summer of 1871, I came back to G\"ottingen with Stolz. [...] He was above all a logician, and during my endless debates with him, the idea that the non-Euclidean geometries were part of Cayley's projective geometry became very clear to me. I imposed it on my friend after a stubborn resistance. I formulated this idea in a short note that appeared in the G\"ottingen Nachrichten, and then in a first memoir, which appeared in Volume 4 of the \emph{Annalen}.
\end{quote}
 A couple of pages later, Klein, talking about his second paper \cite{Klein-Ueber1}, writes (see also \cite{Birkhoff-Bennett}):
  \begin{quote}\small
 I investigated in that paper the foundations of von Staudt's [geometric] system, and had a first contact with modern axiomatics. [...] However, even this extended presentation did not lead to a general clarification. [...] Cayley himself mistrusted my reasoning, believing that a ``vicious circle" was buried in it.
\end{quote} 
 Cayley was more interested in the foundational aspect of projective geometry and his approach was more abstract than that of Klein. In a commentary on his paper \cite{Cayley1859} in his \emph{Collected mathematical papers} edition \cite{Cayley-collected} (Vol. II, p. 605), he writes:
  \begin{quote}\small
  As to my memoir, the point of view was that I regarded ``coordinates" not as distances or ratios of distances, but as an assumed fundamental notion not requiring or admitting of explanation. It recently occurred to me that they might be regarded as mere numerical values, attached arbitrarily to the point, in such wise that for any given point the ratio $x:y$ has a determinate numerical value, and that to any given numerical value of $x:y$ there corresponds a single point. And I was led to interpret Klein's formul\ae \ in like manner; viz. considering $A,B,P, Q$ as points arbitrarily connected with determinate numerical values $a,b,p,q$, then the logarithm of the formula would be that of $(a-p)(b-q) \div (a-q)(b-q)$. But Prof. Klein called my attention to a reference (p. 132 of his second paper) to the theory developed in Staudt's Geometrie der Lage, 1847. The logarithm of the formula is $\log (A,B,P,Q)$ and, according to Staudt's theory $(A,B,P,Q)$, the anharmonic ratio of any four points, has independently of any notion of distance the fundamental properties of a numerical magnitude, viz. any two such ratios have a sum and also a product, such sum and product being each of them like a ratio of four points determinable by purely descriptive constructions.  \end{quote}
 Cayley refers here to von Staudt's notion of a point as a \emph{harmonic conjugate} relatively to three  other points, a definition which was also meant to be independent of any notion of distance (\cite{Staudt} p. 43).
  
Let us end this section by quoting J. E. Littlewood from his  Miscellany \cite{Littlewood}, where he stresses the importance of Cayley's idea: 
 \begin{quote} \small The question recently arose in a conversation whether a dissertation of 2 lines could deserve and get a Fellowship. 
I had answered this for myself long before; in mathematics 
the answer is yes.  Cayley's projective definition of length is a clear case if we may interpret ``2 lines" with reasonable latitude. With Picard's Theorem\footnote{Littlewood is talking here about Picard's theorem saying that if $f:\mathbb{C}\to \mathbb{C}$ is an entire and non-constant function, then it is either surjective or it misses only one point.} it could be literally 2, one of statement, 
one of proof. [...]  With Cayley the importance of the idea is obvious at first sight."\footnote{Littlewood adds: With Picard the situation is clear enough today 
(innumerable papers have resulted from it).  But I can 
imagine a referee's report: ``Exceedingly striking and a 
most original idea. But, brilliant as it undoubtedly is, it 
seems more odd than important; an isolated result, unrelated to anything else, and not likely to lead anywhere."}
\end{quote}

Finally, we point out to the reader that when Cayley talks about a metric space, he does not necessarily mean a metric space as we intend it today. We recall that the axioms of a distance in the sense of a metric, as we intend them today, were formulated by Maurice Fr\'echet (1878-1973)\index{Fr\'echet (Maurice Ren\'e)} in his thesis, defended in 1906. The idea of a ``metric" was somehow vague for Cayley and Klein. 

\section{Beltrami and the Beltrami-Cayley-Klein\index{Beltrami-Cayley-Klein model} model of the hyperbolic plane}
Beltrami's Euclidean model for hyperbolic geometry was a major element in the development of that theory.  Although he did not write any major text on the relation between non-Euclidean and projective geometry, Beltrami was well aware of the works of Cayley and Klein, and, in fact, he was not far from being one of the main actors in this episode. Let us start by quoting Klein on Beltrami's involvement in this intricate story.
This is extracted from the introduction to \cite{Klein-Ueber} (Stillwell's translation p. 73):
\begin{quote}\small
Now since it will be shown that the general Cayley measure\index{Cayley measure} in space of three dimensions covers precisely the hyperbolic, elliptic and parabolic geometries, and thus coincides with the assumption of constant curvature, one is led to the conjecture that the general Cayley measure agrees with the assumption of constant curvature in any number of dimensions. This in fact is the case, though we shall not show it here. It allows one to use formulae, in any spaces of constant curvature, which are presented here assuming two or three dimensions. It includes the facts that, in such spaces, geodesics can be represented by linear equations, like straight lines, and that the elements at infinity form a surface of second degree, etc. These results have already been proved by Beltrami, proceeding from other considerations; in fact, it is only a short step from the formulae of Beltrami to  those of Cayley.
\end{quote}
 
In fact, Beltrami,\index{Beltrami (Eugenio)} two years before Klein wrote published his first paper \cite{Klein-Ueber}, wrote the following to Ho\"uel\index{Ho\"uel (Guillaume-Jules)} (letter dated July 29, 1869 \cite{Beltrami-Boi} p. 96-97):
   \begin{quote}\small
   The second thing [I will add] will be the most important, if I succeed in giving it a concrete form, because up to now it only exists in my head in the state of a vague conception, although without any doubt it is based on the truth. This is the conjecture of a straight analogy, and may be an identity, between pseudo-spherical\index{geometry!pseudo-spherical} geometry\footnote{This is the term used by Beltrami to denote hyperbolic geometry.} and the theory of Mr. Cayley on the \emph{analytical origin of metric ratios}, using the \emph{absolute} conic (or quadric). However, since the theory of invariants plays there a rather significant role and because I lost this a few years ago, I want to do it again after some preliminary studies, before I address this comparison.
      \end{quote}
Three years later, in a letter to Ho\"uel, written on July 5, 1872, Beltrami regrets the fact that he let Klein outstrip him (\cite{Beltrami-Boi} p. 165):
\begin{quote}\small
The principle which has directed my analysis\footnote{Beltrami refers here to a note \cite{beltrami-Osservatione} which he had just published in \emph{Annali di Matematica}.} is precisely that which Mr. Klein has just developed in his recent memoir\footnote{Beltrami refers here to Klein's paper \cite{Klein-Ueber}.} on non-Euclidean  geometry, for 2-dimensional spaces. In other words, from the analytic point of view, the geometry of spaces of constant curvature is nothing else than Cayley's doctrine of the absolute. I regret very much to have let Mr. Klein supersede me on that point, on which I had already assembled some material, and it was my mistake of not giving enough weight to this matter. Beside, this point of view is not absolutely novel, and it is precisely for that reason that I was not anxious to publish my remark. It is intimately related to an already old relation of Mr. Chasles concerning the angle between two lines regarded as an anharmonic ratio (Geom. sup. art. 181) \cite{Chasles-Traite}) and to a theorem of Mr. Laguerre Verlay\footnote{Edmond Laguerre-Verlay, or, simply, Laguerre\index{Laguerre (Edmond)} (1834-1886) studied at the \'Ecole polytechnique, and after that he became an officer in the army. In 1883, he was appointed professor at Coll\`ege de France and two years later he was elected at the French Academy  of Sciences. Laguerre was a specialist of projective geometry and analysis. His name is connected with orthogonal polynomials (the \emph{Laguerre polynomials}). He is the author of 140 papers and his collected works were edited by Hermite, Poincar\'e and Rouch\'e.} (Nouv. Ann. 1853, Chasles, Rapport sur les progr\`es de la g\'eom\'etrie, p. 313). All that Cayley did is to develop an analytic algorithm and, above all, to show that in the general geometry, the theory of rectilinear distances responds exactly to that of angle distances in ordinary geometry. He also showed how and under what circumstances the Euclidean theory of distance differs from the general theory, and how it can be deduced from it by going to the limit.
\end{quote}
      Finally, we quote a letter that Beltrami wrote in the same year to D'Ovidio\footnote{Enrico D'Ovidio\index{D'Ovidio (Enrico)} (1842-1933) was an Italian geometer who is considered as the founder of the famous Turin geometry school. Like Klein, he worked on the question of deriving the non-Euclidean metric function from concepts of projective geometry,  paving the way for subsequent works of Giuseppe Veronese, Corrado Segre and others. D'Ovidio was known for his outstanding teaching, his excellent books, and his care for students. Guiseppe Peano, Corrado Segre, Guido Castelnuovo and Beppo Levi were among his students.} (December 25, 1872, cited in \cite{Loria} p. 422-423).
\begin{quote}\small
When I learned about the theory of Cayley, I realized that his \emph{absolute} was precisely this \emph{limit} locus which I obtained from the equation $w=0$, or $x=0$, and I understood that the identity of the results was due to the following circumstance, that is, in (the analytic) projective geometry one only admits \emph{a priori} that the linear equations represent lines of shortest distance, so that in this geometry one studies, without realizing it, spaces of constant curvature. I was wrong in not publishing this observation, which has been made later on by Klein, accompanied by many developments of which, for several of them, I had not thought.
\end{quote}

We saw that in the case where the fundamental conic used in Cayley's construction is real, the measure defined on the interior of the conic gives a model of Lobachevsky's geometry. Klein recovered in this way the model which Beltrami had introduced in his paper \cite{Beltrami-Saggio} in which he noticed that the Euclidean straight lines in the unit disc behave like the non-Euclidean geodesics. It was Klein who provided this model with an explicit distance function, namely, the distance defined by the logarithm of the cross ratio, and he also noticed that the circle, in Beltrami's model, can be replaced by an ellipse.

Although this was not his main goal, Klein used this model to discuss the issue of the non-contradiction of hyperbolic geometry. This was also one of Beltrami's achievements in his paper \cite{Beltrami-Saggio}.\footnote{We recall that this model was discovered by Beltrami four years before he discovered his famous pseudo-spherical model.}  We mention that the non-contradiction issue as well as the relative non-contradiction issue (meaning that if one geometry is contradictory, then the other would also be so) among the three geometries was one of the major concerns of Lobachevsky, see e.g. his \emph{Pangeometry} and the comments in the volume \cite{L}. It is also important to recall that while Beltrami's Euclidean model showed that hyperbolic geometry is consistent provided Euclidean geometry is, Klein's work shows that Euclidean, spherical and hyperbolic geometries are consistent provided projective geometry is consistent. 

We also mention that in \S 14 of the paper \cite{Klein-Ueber}, while he computed the curvature of the metric, Klein obtained the expression, in polar coordinates, of the so-called Poincar\'e metric of the disk.

In his lecture notes \cite{Klein-class} (p. 192), Klein writes the following:\footnote{The English translations of our quotes from  \cite{Klein-class} were made by Hubert Goenner.}
\begin{quote}\small [..] it is the merit of
  Beltrami's Saggio, to emphatically have called attention to the fact
  that the geometry on surfaces of constant negative curvature really
  corresponds to non-euclidean hyperbolic geometry.
  \end{quote}
On p. 240, Klein discusses topology, and says that for 2-dimensional spaces of positive curvature, instead of working like Beltrami on the sphere, where  two geodesic lines intersect necessarily intersect in two points, one can work in elliptic space, where geodesics intersect in only one point.\footnote{ It is interesting that in
the 1928 edition of Klein's course on non-euclidean geometry \cite{Klein-Vorlesungen}, the editor Rosemann removed almost all of Klein's
remarks concerning  Beltrami's contributions made in his course of
1989/90 \cite{Klein-class} while Klein was alive. (This remark was made to the authors by Goenner.)}

 \section{The construction of measures}\label{s:measures}
 We now return to Klein's papers.
  The core of the  paper \cite{Klein-Ueber} starts at \S 3, where Klein describes the construction of one-dimensional projective measures, that is, measures on lines and on circles. The one-dimensional case is the basic case because higher-dimensional measures are built upon this case. Klein refers to the one-dimensional case as the \emph{first kind}.
There are two sorts of \emph{measures}\index{Klein's measure for distances}\index{Klein's measure for angles} to be constructed: measures on points and measures on angles. The measure function on points satisfies the usual properties of a distance function\footnote{It is considered that the first formal statement of the axioms of a distance function as we know the today is due to Fr\'echet in his thesis \cite{Frechet} (1906), but the nineteenth-century mathematicians already used this notion, and they ware aware of geometries defined by distance functions.} except that it can take complex values. The measure for angles is, as expected, defined only up to the addition of multiples of $2\pi$, and at each point it consists of a measure on the pencil of lines that pass through that point. It can also take complex values. Klein specifies the following two properties that ought to be satisfied by measures: 
\begin{enumerate}
\item the measures for points satisfy an additivity property for triples of points which are aligned; 
\item the two measures (for points and for lines) satisfy the property that \emph{they are not altered by a motion in space}. 
\end{enumerate}
Property (a) says that the projective lines are geodesics for these measures. A metric that satisfies this property is called (in modern terms) \emph{projective}. The motions that are considered in Property (b) are the projective transformations that preserve a conic, which is termed the \emph{basic figure}. Klein then addresses the question of the classification of measures, and he notes that this depends on the classification of the transformations of the basic figures, which in turn depends  on the number of fixed points of the transformation. Since the search for fixed points of such transformations amounts to the search for solutions of a degree-two equation, the transformations that preserve the basic figure fall into two categories:  
\begin{enumerate}
\item Those that fix two (real or imaginary) points of the basic figure, and this is the generic case. They are termed measures \emph{of the first kind}.
\item Those that fix one point of the basic figure. They are termed measures \emph{of the second kind}.
\end{enumerate}

  Klein describes in detail the construction of measures on lines. The overall construction amounts to a division of the circle (seen as the projective line) into smaller and smaller equal parts, using a projective transformation. Thus, if we set the total length of the circle to be 1, the first step will provide two points at mutual distance $\frac{1}{2}$, the second step will provide three points at mutual distance $\frac{1}{3}$, and so forth. Passing to the limit, we get a measure on the circle which is invariant by the action of the given projective transformation.

More precisely, Klein starts with a transformation of the projective line of the form $z\mapsto \lambda z$, with $\lambda$ real and positive. The transformation has two fixed points, called \emph{fundamental elements}, the points $0$ and $\infty$. Applying the transformation to a point $z_1$ on the line, we obtain the sequence of points $z_1,\lambda z_1, \lambda^2 z_1, \lambda^3 z_1, \ldots$. In order to define the measure, Klein divides, for any integer $n$, the line into $n$ equal parts using the transformation $z'=\lambda^{\frac{1}{n}}z$. The $n$th root determination is chosen in such a way that $\lambda^{\frac{1}{n}}z$ lies between $z$ and $\lambda z$. The distance between two successive points is then defined as the $\frac{1}{n}$th of the total length of the line. Iterating this construction, for any two integers $\alpha$ and $\beta$, the distance between $z_1$ and a point of the form $\lambda^{\alpha+\frac{\beta}{n}}z_1$ is set to be the exponent $\alpha+\frac{\beta}{n}$, that is, the logarithm of the quotient $\frac{\lambda^{\alpha+\frac{\beta}{n}}z_1}{z_1}$ divided by $\log \lambda$. By continuity, we can then define the distance between two arbitrary points $z$ and $z_1$ to be the logarithm of the quotient $\frac{z}{z_1}$ divided by the constant $\log \lambda$. The constant $\frac{1}{\log \lambda}$ is denoted henceforth by $c$. 

Klein shows that the measure that is defined in this way is additive, that the distance from a point to itself is zero, and that the distance between two points is invariant by any linear transformation that fixes the fundamental elements $z=0$ and $z=\infty$.
He then observes that the quotient  $\frac{z}{z'}$ may be interpreted as the cross ratio of the quadruple $0,z,z',\infty$. 
Thus, the distance between two points $z,z'$ is a constant multiple of the logarithm of the cross-ratio of the quadruple $0,z,z',\infty$.  In particular, the distance between the two fundamental elements is infinite.
   
In \S 4 of his paper, Klein extends the distance function  $c\log\frac{z}{z'}$ to pairs of points on the complex line 
joining the points $0$ and $\infty$, after choosing a determination of the complex logarithm. He then gives an expression for a general result where he assumes, instead of the special case where the two fundamental elements are $0$ and $\infty$, that these points are the solutions of a second-degree equation
\[\Omega= az^2+2bz+c=0.\]
For two arbitrary points given in homogeneous coordinates, $(x_1,x_2)$ and $(y_1,y_2)$, setting
\[\Omega_{xx}= ax_1^2+2bx_1x_2+cx_2^2,\]
\[\Omega_{yy}= ay_1^2+2by_1y_2+cy_2^2\]
and
\[\Omega_{xy}= ax_1y_1+2b(x_1y_2+x_2y_1)+cx_2y_2,\]
the distance between the two points is 
\begin{equation}\label{eq:Omega0}
 c \log \frac{\Omega_{xy}+\sqrt{\Omega_{xy}^2-\Omega_{xx}\Omega_{yy}}}{\Omega_{xy}-\sqrt{\Omega_{xy}^2-\Omega_{xx}\Omega_{yy}}}.
\end{equation}

Later on in Klein's paper, the same formula, with the appropriate definition for the variables, defines a measure between angles between lines in a plane and between  planes in three-space.

In \S 5, Klein derives further properties of the construction of the measure, distinguishing the cases where the two fundamental points are respectively real distinct, or conjugate imaginary, or coincident. 

The first case was already treated; the two fundamental points are infinite distance apart, and they are both considered at infinity. He observes that this occurs in hyperbolic geometry, a geometry where any line has two points at infinity. 

The case where the two fundamental points are conjugate imaginary occurs in elliptic geometry where a line has no point at infinity. Klein shows that in this case all the lines are finite and have a common length, whose value depends on the constant $c$ that we started with.  The distance \ref{eq:Omega0} between two points becomes
\begin{equation}\label{eq:Omega}
 2 i c \arccos\frac{\sqrt{\Omega_{xy}}}{{\sqrt{\Omega_{xx}\Omega_{yy}}}}.
\end{equation}
Klein notes that a particular case of this formula appears in Cayley's paper, who used only the value $-\frac{i}{2}$ for $c$ and when, consequently, the term in front of $\arccos$ is equal to 1.

Klein studies the case where the two fundamental points coincide in \S 6. This case concerns Euclidean (parabolic) geometry.\index{geometry!parabolic}\index{geometry!Euclidean} It is more complicated to handle than the other cases and it needs a special treatment. One complication arises from the fact that in this case Equation (\ref{eq:Omega}), which has a unique solution, leads to distance zero between the points $x$ and $y$. The problem is resolved by considering this case as a limit of the case where the equation has two distinct solutions. Klein derives from there the formula for the distance on a line in which there is a unique point at infinity, that is, a unique point which is infinitely far from all the others.

In \S 7, Klein introduces a notion of \emph{tangency of measures at an element}. For this, he introduces two measures associated to a basic figure of the first kind, which he calls ``general" and ``special", and which he terms as ``tangential". The overall construction amounts to the definition of infinitesimal geometric data, and it is also used to define a notion of a \emph{curvature of a general measure}. The sign and the value of this curvature depend on some notion of deviation, which he calls ``staying behind or running ahead" of the general measure relative to the special measure. He shows that the value of this geometrically defined curvature is constant at every point, and equal to $\-\frac{1}{4c^2}$, where $c$ is the characteristic constant of the general measure. Using Taylor expansions, Klein shows that the three geometries (elliptic, parabolic and hyperbolic) are tangentially related to each other, which is a way of saying that infinitesimally, hyperbolic and spherical geometry are Euclidean. The value of $c$ is either real or imaginary so that one can get positive or negative curvature.\footnote{The result should be real, and for that reason, $c$ has sometimes to be taken imaginary. This is to be compared with the fact that some (real) trigonometric functions can be expressed as functions with imaginary arguments.}

 In \S 8, Klein outlines the construction of the measure for basic figures of the second kind, that is, measures on planes and measures on pencils of dihedral angles between planes. He uses for this an auxiliary conic. This is the so-called \emph{fundamental conic} (the conic that is called the \emph{absolute} by Cayley). Each projective line intersects this conic in two points (real, imaginary or coincident). The two points play the role of fundamental points for the determination of the metric on that line, and the problem of finding a measure is reduced to the 1-dimensional case which was treated before. The fundamental conic is the locus of points which are   infinitely distant from all others. 
 
Measures on rays in the plane are based on the fact that at each point, there are rays that start at that point and that are are tangent to the conic. Again, these rays are solutions of a certain quadratic equation and they may be distinct real, distinct imaginary or coincident. The two tangent rays are taken to be the fundamental rays for the angle determination in the sense that the angle between two arbitrary rays is then taken to be the cross ratio of the quadruple formed by these rays and the fundamental rays. The multiplicative constant is not necessarily the same in the formulae giving the measures on lines and on rays.

 Klein then determines an analytic expressions for these measures. It turns out that the formulae are the same as those obtained in \S 4. If the equation of the fundamental conic is 
 \[\Omega=\sum_{i,j=1}^3 a_{ij}x_i y_j=0,\]
  then
  the distance between the two points $x$ and $y$, in homogeneous coordinates $(x_1,x_2,x_3)$ and $(y_1,y_2,y_3)$, is  
   \begin{equation}\label{eq:Omega1}
 c \log \frac{\Omega_{xy}+\sqrt{\Omega_{xy}^2-\Omega_{xx}\Omega_{yy}}}{\Omega_{xy}-\sqrt{\Omega_{xy}^2-\Omega_{xx}\Omega_{yy}}}
\end{equation}
where $\Omega_{xx}, \Omega_{yy}$, etc. are the expressions obtained by substituting in $\Omega$ the coordinates $(x_1,x_2, x_3)$ of a point $x$ or $(y_1,y_2, y_3)$ of a point $y$, etc.
Equivalently,  we have
\begin{equation}\label{eq:Omega2}
2 i c \arccos\frac{\sqrt{\Omega_{xy}}}{{\sqrt{\Omega_{xx}\Omega_{yy}}}}.
\end{equation}
That is, one obtains again Formulae (\ref{eq:Omega0}) and  (\ref{eq:Omega}) of \S 4.

Concerning measures on angles, the equations have a similar form.
One takes the equation of the fundamental conic in line coordinates to be
 \[\Phi_{u,v}=\sum_{i,j=1}^3 A^{ij}u_i v_j=0.\]
The distance between the two points $u$ and $v$ in homogeneous coordinates $(u_1,u_2,u_3)$ and $(v_1,v_2,v_3)$  is then
   \begin{equation}\label{eq:Phi1}
 c' \log \frac{\Phi_{uv}+\sqrt{\Phi_{uv}^2-\Phi_{uu}\Phi_{vv}}}{\Phi_{uv}-\sqrt{\Phi_{uv}^2-\Phi_{uu}\Phi_{vv}}}
\end{equation}
or, equivalently, 
\begin{equation}\label{eq:Phi2}
2 i c' \arccos\frac{\sqrt{\Phi_{uv}}}{{\sqrt{\Phi_{uu}\Omega_{vv}}}}.
\end{equation}
where $\Omega_{uu}, \Omega_{vv}$ have the same meaning as before.

The constant $c'$ is  in general different from $c$. In general, the constants are chosen so that the result is real. 

The measures on points and on lines are defined by similar formulae. This is a consequence of the fact that they are solutions of second-degree equations, and that the coefficients of the two equations are related to each other by the duality in projective geometry. Duality is discussed in the next section. 

 \S 9 concerns the properties of the projective transformations of the plane that preserve a conic.  Klein points out that there is a ``threefold infinity" of such transformations (in other words, they form a 3-dimensional group), and he starts a classification of such transformations, based on the fact that each transformation fixes two points of the conic and reasoning on the line connecting them, on the tangents at these points, on their point of intersection, and working in the coordinates associated to the triangle formed by the connecting line and the two tangents. The classification involves the distinction between real conics with real points and real conics without real points. The aim of the analysis is to prove that the transformations that map the conic into itself preserve the metric relations between points and between angles. There is also a polar duality determined by the conic. With this duality, a quadruple formed by two points and the intersection of the line that joins them with the conic corresponds to a quadruple formed by two lines and tangents to the conic that pass through the same point. This correspondence preserves cross ratios. The duality is such that the distance between two points is equal to the angle between the dual lines. This is a generalization of the polar duality that occurs in spherical geometry.

After the discussion of projective measures  between points in \S 9,  Klein considers in \S 10 measures for angles in pencils of lines and of planes. In this setting, he uses a \emph{fundamental cone of second degree} instead of the fundamental conic. He also appeals to polarity theory and he makes a relation with the measure obtained in the previous section. The result brought out at the end of the previous section is interpreted here as saying that the angle between two planes is the same as the angle between their normals, and this has again an interpretation in terms of spherical geometry duality.

In \S 11, Klein develops a model for spherical geometry\index{geometry!spherical} that arises from his measures associated to conics.

 When the fundamental conic is imaginary, setting $c=c_1\sqrt{-1}$ and $c'=c'_1\sqrt{-1}$, the measures for lines and for angles\index{Klein's measure for distances}\index{Klein's measure for angles} are found to be respectively 
\begin{equation} \label{FK1}
2c_1\arccos \frac{xx'+yy'+zz'}{\sqrt{x^2+y^2+z^2} \sqrt{x'^2+y'^2+z'^2}}
\end{equation}
and
\begin{equation} \label{FK2}
2c'_1\arccos \frac{uu'+vv'+ww'}{\sqrt{u^2+v^2+w^2} \sqrt{u'^2+v'^2+w'^2}},
\end{equation}
which are the familiar formulae for angle measure on a sphere.
In particular,  the distance between any two points is bounded, as expected. In fact,all lines are closed, they have finite length, and these lengths have a common value, $2c_1\pi$, which is (up to a constant multiple) the angle sum of a pencil, which is $2c'_1\pi$. The point measure is completely similar to the angle measure. This again can be explained the duality  between points and lines in spherical geometry. Klein concludes from this fact that ``plane trigonometry, under this measure, is the same as spherical trigonometry" and that ``the plane measure just described is precisely that for elliptic geometry". By choosing appropriately the constants $c_1$ and $c'_1$, the angle sum of any pencil becomes $\pi$ and the maximal measure between points becomes also $\pi$. Klein also deduces that in that geometry, the angle sum of a plane triangle is greater than $\pi$, as for spherical triangles, and only equal to $\pi$ for infinitesimally small triangles.

In \S 12, Klein describes the construction that leads to hyperbolic geometry. This is the case where the absolute is a real fundamental conic in the plane. In this case, the constant $c$ that appears in the general formula \ref{eq:Omega0} for distances is taken to be real. The points in the plane are divided into three classes: the points inside the conic, the points on the conic and the points outside the conic. The points inside the conic are those that admit no real tangent line to the conic.
The points on the conic are those that admit one real tangent line. The points points outside the conic are those that admit two real tangent lines.

Likewise, the lines in the plane are divided into three classes: the lines that meet the conic in two real points, those that meet the conic in a unique (double) real point and those that do not meet the conic in any  real point. Klein claims that this case corresponds to hyperbolic geometry. To support this claim, he writes:

\begin{quote}\small
The geometry based on this measure \emph{corresponds completely with the idea of hyperbolic geometry}, when we set the so far undetermined contant $c'_1$ equal to $\frac{1}{2}$, making the angle sum of a line pencil equal to $\pi$. In order to be convinced of this, we consider a few propositions of hyperbolic geometry in somewhat more detail.
\end{quote}
The propositions that Klein considers are the following:
\begin{enumerate}
\item Through a point in the plane there are two parallels to a given line, i.e. lines meeting the points at infinity of the given line.
\item The angle between the two parallels to a given line through a given point decreases with the distance of the point from the line, and as the point tends to infinity, this angle tends to $0$, i.e. the angle between the two parallels tends to zero.
\item The angle sum of a triangle is less than $\pi$. For a triangle with vertices at infinity, the angle sum is zero.
\item Two perpendiculars to the same line do not meet.
\item A circle of infinite radius is not a line.
\end{enumerate}

Klein notes that these properties are satisfied by his geometry. This is not a full proof of the fact that the geometry defined using the distance function he described \emph{is} hyperbolic geometry, but it is a strong indication of this fact. In fact, it is surely possible, but very tedious, to show that all the axioms of  hyperbolic geometry are satisfied by his geometry. 
Klein then adds (Stillwell's translation p. 99):
\begin{quote}\small
Finally, the \emph{trigonometric formulae} for the present measure are obtained immediately from the following considerations. In \S 11 we have seen that, on the basis of an imaginary conic in the plane and the choice of constants $c=c_1i$, $c'=c'_1i=\frac{\sqrt{-1}}{2}$, the trigonometry of the plane has the same formulae as spherical trigonometry when one replaces the sides by sides divided by $2c_1$. The same still holds on the basis of a real conic. Because the validity of the formulae of spherical trigonometry rests on analytic identities that are independent of the nature of the fundamental conic. The only difference from the earlier case is that $c_1=\frac{c}{i}$ is now imaginary.
\\
\emph{The trigonometric formulae that hold for our measure result from the formulae of spherical trigonometry by replacing sides by sides divided by $\frac{c}{i}$.}
\\
But this is the same rule one has for the trigonometric formulae of hyperbolic geometry. The constant $c$ is the characteristic constant of hyperbolic geometry. One can say that planimetry, under the assumption of hyperbolic geometry, is the same as geometry of a sphere with the imaginary radius $\frac{c}{i}$.
\\The preceding immediately gives a model of hyperbolic geometry, in which we take an arbitrary real conic and construct a projective measure on it. Conversely, if the measure given to us is representative of hyperbolic geometry, then the infinitely distant points of the plane form a real conic enclosing us, and the hyperbolic geometry is none other than the projective measure based on this conic.
\end{quote}

\S 13 concerns parabolic geometry.\index{geometry!parabolic} In this case, the fundamental figure at infinity is a degenerate conic. It is reduced to a pair of points, and, in Klein's words, it constitutes a ``bridge between a real and an imaginary conic section". The metric obtained is that of Euclidean geometry and  line joining the pair of points at infinity (the degenerate conic) is the familiar line at infinity of projective geometry. In this sense, parabolic geometry is regarded as a transitional geometry, sitting between hyperbolic and elliptic geometry. To understand how this occurs, Klein gives the example of a degeneration of a hyperbola. A hyperbola has a major and a minor axis, which are symmetry axes, the major axis being the segment joining the two vertices $a$ and $-a$ (and the length of this axis is therefore equal to $2a$) and the minor axis being perpendicular to the major one, with vertices at points $b$ and $-b$,  of length $2b$. The major and minor axes are also the two perpendicular bisectors of the sides of a rectangle whose vertices are on the asymptotes of the hyperbola. 
In a coordinate system where the two axes are taken as the major and major axes, the equation of the hyperbola is $\frac{x^2}{a^2}-\frac{y^2}{b^2}= 1$. The minor axis is also called the imaginary axis because of the minus sign occurring in this equation. 

The degeneration of the hyperbola into two imaginary points is obtained by keeping fixed the imaginary axis and shrinking to zero the major axis. Meanwhile, the two branches of the hyperbola collapse to the line carried by the minor axis, covering it twice. This line represents a degenerate conic, and in fact, as Klein points out, insofar as it is enveloped by lines, it is represented by the two conjugate imaginary points. The associated measure on the plane is called a \emph{special measure}, because it uses a pair of points instead of a fundamental conic. Klein obtains an analytic formula that gives the associated distance between points.
Starting with the general expression 
\[2ic \arcsin \frac{\sqrt{\Omega_{xy}^2-\Omega_{xx}\Omega_{yy}}}{\sqrt{\Omega_{xx}\Omega_{yy}}},\]
where $\Omega=0$ is the equation of the conic and $\Omega_{xx}$, $\Omega_{xy}$ and $\Omega_{yy}$ are as he defined in \S 4, and taking limits when the conic degenerates to a pair of points, he deduces that the distance between two points $(x,y)$ and $(x',y')$ is 
\[\frac{C}{k^2}\sqrt{(x-x')^2+(y-y')^2}\]
which up to a constant factor is the Euclidean distance in coordinates. He concludes this section with the following: 
\begin{quote} \small
We want to stress that with imaginary fundamental points the trigonometric formulae become the relevant formulae of parabolic geometry, so the angle sum of a triangle is exactly $\pi$, whereas with a real fundamental conic it is smaller, and with an imaginary conic it is larger.
\end{quote}

In \S 14, Klein considers again the notion of a ``measure on a plane which is tangent to a general measure at a point", that is, he considers infinitesimal distances. This leads him again to the definition of a notion of curvature which is equivalent to Gaussian curvature. Klein then uses a duality, where the dual of a point of the given geometry is a ``line at infinity" which is the polar dual of the given point with respect to the fundamental conic. He obtains a qualitative definition of curvature which turns out to be equivalent to the Gaussian curvature. This leads to the conclusion that the curvature of a general measure is the same at all points and is equal to $\frac{-1}{4c^2}$, that it is positive if the fundamental conic is imaginary (the case of elliptic geometry) and negative if the fundamental conic is real (the case of hyperbolic geometry). In the transitional case (parabolic geometry), which is a limiting case in which the fundamental conic degenerates into a pair of imaginary points,  the curvature is zero.
Klein concludes this section with the following statement: 

\emph{According to whether we adopt the hypothesis of an elliptic, hyperbolic or parabolic geometry, the plane is a surface with constant positive, constant negative or zero curvature}.

In \S 15, Klein talks about a continuous transition from hyperbolic to parabolic and from spherical to parabolic. Let us quote him (Stillwell's translation p. 107):
\begin{quote}\small
If we are actually given parabolic geometry we can immediately constructs a geometry which models hyperbolic geometry by constructing a general measure with real fundamental conic, tangential to the given special measure at a point of our choice. We achieve this by describing a circle of radius $2c$ centred on our point, and using it as the basis for a projective measure with the constant $c$ determining the distance between two points and the constant $c'=\frac{\sqrt{-1}}{2}$ determining the angles between two lines. This general measure approaches the given parabolic measure more closely as $c$ becomes larger, coinciding with it completely when $c$ becomes infinite.
\\
In a similar way we construct a geometry that shows how elliptic geometry can tend toward the parabolic. To do this it suffices to give a pure imaginary value $c_1i$ to the $c$ we used previously. Then we fix a point at distance $2c_1$ above the given point of contact and take the distance between two points of the plane to be $c_1$ times the angle the two points subtend at the fixed point. The angle between the two lines in the plane is just the angle they subtend at the fixed point. The resulting measure approaches more closely to the parabolic measure the greater $c_1$ is, and it becomes equal to it when $c_1$ is infinite. 
\\
When elliptic or hyperbolic geometry is actually the given geometry one can in this way make a model presenting its relationship with parabolic or the other geometry.
\end{quote}

\S 16 concerns projective measures in space. The same procedure as before is used, with a second-degree  fundamental surface in 3-space instead of the fundamental curve in the plane. The case where the fundamental surface is imaginary leads to elliptic geometry. The case where it is \emph{real} and \emph{not ruled}, leads to hyperbolic geometry.\index{geometry!hyperbolic}  The case where the fundamental surface degenerates to a conic section leads to parabolic geometry, and this conic section becomes the \emph{imaginary circle at infinity}. The case where the fundamental surface is \emph{real} and \emph{ruled}, that is, a one-sheeted hyperboloid, is not related to any of the classical geometries, it leads to a geometry which is not locally Euclidean but pseudo-Euclidean.

 The title of  \S 17 is ``The independence of projective geometry from the theory of parallels". Klein observes that projective geometry insofar as it uses the notions of homogeneous coordinates and the cross ratio, is defined in the setting of parabolic geometry. He notes however that in the same way as one can construct projective geometry starting with parabolic geometry, one can construct it also on the basis of hyperbolic and elliptic geometry. He then notes that projective geometry can be developed without the use of any measure, using the so-called incidence relations, referring to the work of von Staudt.
 
In the conclusion to the paper (\S 18), Klein notes that by a consideration of the sphere tangent to the fundamental surface, one is led to only the three geometries considered, elliptic, hyperbolic and the transitional one, that is, the parabolic.

     \section{Klein's second paper} \label{s:second-p}
  The second paper (\cite{Klein-Ueber1} same title, Part II), appeared two years after the first one. It has a more general character, it is in the spirit of his \emph{Erlangen program}, and it is less technical than the first paper. There does not seem to be an available English translation of that paper. In this paper, Klein gives some more details on results he obtained in the first paper. Let us quote Klein from his G\"ottingen lecture notes of Klein \cite{Klein-class} (p. 286-287; Goenner's translation):
  \begin{quote}\small
  When the accord of Cayley's measure
  geometry and non-euclidean geometry was stated *, it became
  essential to draw conclusions from it. On these conclusions I wish
  to attach most importance, although they were developed more in
  detail only in the second paper, when  I noticed that the very same
  [conclusions] did appear to other mathematicians not as self-evident
  as for myself. [..] (footnote) *  Beltrami and Fiedler also had
  noticed this accord, as they later wrote to me. 
  \end{quote}

  We now give a brief summary of the content of the paper \cite{Klein-Ueber1}.  
  
  This paper has two parts. In the first part, Klein develops the idea of a transformation group that characterizes a geometry. In the second part, he develops an idea concerning projective geometry which he had also mentioned in the first paper \cite{Klein-Ueber}, namely, that projective geometry is independent from the Euclidean parallel postulate (and from Euclidean geometry).\footnote{One of the basic features of projective geometry is that in the arguments that involves lines, unlike in Euclidean geometry, one does not have to distinguish between the cases where the lines intersect or are parallel. In projective geometry, any two distinct points define one line, and any two distinct lines intersect in one point. We already mentioned that this principle is at the basis of duality theory in projective geometry.}  Klein insisted on this fact, because, as he wrote, some mathematicians thought that there was a vicious circle in his construction of Euclidean geometry from projective geometry, considering that the definition of the cross ratio uses Euclidean geometry, since it involves a compounded ratio between four Euclidean segments.  Some also thought that there was a contradiction in Klein's reasoning, since in spherical and hyperbolic geometry Euclid's parallel axiom is not satisfied, so a formula for the metric defining these geometries cannot be based on the distance function of Euclidean geometry where the parallel axiom is satisfied.  In fact, as was already recalled above, in his \emph{Geometrie der Lage} \cite{Staudt}, von Staudt had already worked out a purely projective notion of the cross ratio, independent of any notion of distance.\footnote{We can quote here Klein from his \emph{Erlangen program} \cite{Klein-Erlangen}: ``We might here make mention further of the way in which \emph{von Staudt} in his \emph{Geometrie der Lage} (N\"urnberg, 1847) develops projective geometry, -- i.e., that projective geometry which is based on the group containing all the real projective and dualistic transformations". And in a note, he adds: ``The extended horizon, which includes \emph{imaginary} transformations, was first used by \emph{von Staudt} as the basis of his investigation in his later work, \emph{Beitr\"age zur Geometrie der Lage} (N\"urnberg, 1856-60)." For the work of von Staudt and in his influence of Klein, and for short summaries of the  \emph{Geometrie der Lage}  and the \emph{Beitr\"age zur Geometrie der Lage}, we refer the reader to the paper \cite{Gray-K} in this volume.}
   In his \emph{Lectures on the development of mathematics in the XIXth century} (1926-1927),  \cite{Klein19} Klein returns to the history and he writes the following (Vol. 1, p. 153):
  \begin{quote}\small
  More important is the objection I received from mathematicians. In my paper written in Volume IV of the \emph{Annalen}, I did not expect the logical difficulties that the problem raised, and I had started an innocent use of metric geometry. It is only at the end that I mentioned in a very brief way the independence of projective geometry from any metric, referring to von Staudt. I was accused from everywhere of making circular reasoning. The purely projective definition of von Staudt of the cross ratio as a number was not understood, and people stood firmly on the idea that this number was only given as a cross ratio of four Euclidean numbers.
   \end{quote}
   
   We also quote Cayley's citation of R. S. Ball \cite{Cayley-collected} ( Vol. II, p. 605):
  
  \begin{quote}\small
  I may refer also to the memoir, Sir R. S. Ball ``On the theory of content," \emph{Trans. R. Irish Acad.} vol. {\sc xxix} (1889), pp. 123--182, where the same difficulty is discussed. The opening sentences are -- ``In that theory [Non-Euclidian geometry] it seems as if we try to replace our ordinary notion of distance between two points by the logarithm of a certain anharmonic ratio. But this ratio itself involves the notion of distance measured in the ordinary way. How then can we supersede the old notion of distance by the non-Euclidian notion, inasmuch as the very definition of the latter involves the former?
\end{quote}

From this, let us conclude two different things:
\begin{enumerate}
\item There was a great deal of confusion about Klein's ideas, even among the most brilliant mathematicians.
\item The mathematicians were not only interested in formulae, but they were digging in the profound meaning that these formulae express.
\end{enumerate}

  We end this paper with a brief summary of the content of \cite{Klein-Ueber1}, since no available translation exists. The reader can compare the content of this paper with the summary of the \emph{Erlangen program} lecture given in Chapter 1 of this volume  \cite{Gray-K}.

 The introduction contains historical recollections on the works of Cayley and von Staudt which, according to Klein, did not have yet any applications. Klein then recalls that the ``problem of parallels", that is, the problem of deciding whether Euclid's parallel axiom follows or not from the other axioms of Euclidean geometry was settled. He mentions the works of the founders of modern geometry, and he says that each of them brought new mathematical concepts, in particular, new examples of spaces of constant curvature. At the same time, several open questions remain to be solved, and other things need to be made more precise. Klein then mentions the spaces of variable curvature constructed by Riemann. All these works contribute to new points of view on spaces and on mechanics.  He also recalls that there is a difference between the metrical and the projective points of view, and  he declares that the geometries of constant curvature should be simpler to study.

In Section 1 of the first part of the paper, Klein considers the concept of higher dimensions. He  mentions the relation between constant curvature manifolds and projective manifolds.\footnote{In \cite{Klein-Erlangen}, the term ``Mannigfaltigkeiten", which for simplicity we translate by ``manifold", is usually  translated by ``manifoldness". See the comments in \cite{Gray-K} in this volume on the meaning of the word manifoldness.}  Analytic geometry allows the passage to higher dimensions, working in analogy with the low dimensions that we can visualize. He points out that on a  given line, we can consider either the real points or all points. He then recalls the definition of the cross ratio. 

Section 2 concerns transformations. Klein explains the notion of composition of transformations, and he considers in particular the case of collinea\-tions, forming a group. He then presents the idea of group isomorphism. The reader should recall that these ideas were relatively new at that time. 

 Section 3 concerns ``invariant", or ``geometric", properties. A property is geometric if it is independent of the location in space. A figure should be indistinguishable from its symmetric images. The properties that we seek are those that remain invariant by the transformations of the geometry.
 
 In Section 4, Klein develops the idea that the methods of a given geometry are characterized by the corresponding groups. This is again one of the major ideas that he had expressed in his \emph{Erlangen program}. He elaborates on the significance of projective geometry, and in particular on the transformations that leave invariant the imaginary circle at infinity. The methods depend on the chosen transformation group.

The discussion on   Sections 3 and is also confirmed by Klein in Klein's lecture notes of 1898/90 (\cite{Klein-class}, p. 120; Goenner's translation):
 \begin{quote}\small 
 In  contradistinction [to Helmholtz], I had the generic thought that, in studying manifolds under the viewpoint of giving them a {\em geometric} character, one can put ahead any transformation group  [..].* ~Above all it is advisable to chose the collineations (linear
  transformations) as such a group. [..] This then is the specially
  so-called {\em invariant theory}. - (footnote) * Ann. VI, p. 116 et
  seq., as well as notably the Erlangen program. 
  \end{quote}  
  
 Section 5 concerns generalizations to higher-dimensional spaces. The simplest transformation groups are the groups of linear transformations. They give rise to projective geometry. Although there is no distance involved, this is considered as a geometry.  Klein introduces the word ``invariant theory", where we have no distance involved, but we look for invariant objects. Modern algebra is helpful in that study. In the case where we have a metric, we have an invariant quadratic form. He declares that he will study the case where there is none. He introduces a notion of differential of a map. At the infinitesimal level, the differential behaves like a linear map.
 
 In Section 6, Klein considers spaces of constant nonzero curvature. He refers to Beltrami, who showed that in such a space we can define geodesics by equations that are linear in the appropriate  coordinates. He raises the question of understanding the transformations of a manifold of constant curvature in the projective world, and this is done in linearizing them. Indeed, by choosing adequate coordinates, the group of transformations that we attach to a manifold of constant curvature is contained in the group of linear transformations. (\emph{``The transformation group of a constant curvature geometry is reducible to a transformation group which preserves a quadratic form".}) Klein says however that there is a difference between his viewpoint and the one of Beltrami, namely, Klein starts with complex variables and then he restricts to real variables. This gives a uniform approach to several things. 
 
 Elliptic space is obtained from the sphere by identifying antipodal points so that there is a unique geodesic connecting two points. In higher dimensions similar objects exist. 
 
 Sections  7 and 8 concern the description of constant curvature manifolds in terms of projective notions, and Klein recalls the definition of the distance using the cross ratio.  
  
  In Section 9, Klein defines a point at infinity of the space as a representative of a class of geodesic lines. 
  
The subject of the second part of the paper is the fact that, following von Staudt, one can construct projective space independently of the parallel axiom.
 
 In the first section of this part, Klein explains various constructions that are at the basis of projective geometry. He also introduces  the betweenness relation. He talks about lines and pencils of planes, of the cross ratio and of the notion of harmonic division, and he states the fact that there is a characterization of 2-dimensional projective geometry. He recalls that von Staudt, in his work, used the parallel axiom, but that without essential changes one can recover the bases of projective geometry without the parallel axiom. He then studies the behavior of lines and planes, and the notion of asymptotic geodesics.  This section also contains a detailed discussion of von Staudt's axiomatic approach. 
 
 Section 2  concerns the ``formulation of a proposition which belongs to the general theory of \emph{Analysis situs}". Klein explains how one can attach coordinates to points. 
 
 In Section 3, he returns to the bases of systems of planes and their intersection.
 
  In Section 4, he elaborates on the notion of harmonic element, and on the notion of  betweenness among points. 

 In Section 5, Klein expands on the work of van Staudt on projective transformations.
 
 In Section 6, he talks about pencils of planes and about duality and ideal points;  a point at infinity defines as a class of lines which do not intersect. 
 
 Section 7 concerns the cross ratio and homogeneous coordinates.
 
 In Section 8, Klein gives an analytical proof of the main theorem of projective geometry. He says that there is a characterization of 2-dimensional projective geometry.
 
\medskip
   
To conclude this section, let us insist on the fact that beyond their immediate goal (which is an important  one), the two papers by Klein are full of interesting historical comments and references to works of other mathematicians. They are the expression of the elegant style and the great erudition which characterizes Klein's writings in general.  The reader should remember that in 1871, at the time he wrote the paper \cite{Klein-Ueber}, Klein was only 22 years old.

  \section{Poincar\'e}\label{s:P}

 In his paper \cite{Poincare1887} (1887),\index{Poincar\'e (Henri)} Poincar\'e describes a construction of a set of geometries, using quadrics in three-space.\index{geometry!quadratic surface}\index{quadric} The theory of associating a geometry to a quadric is of course related to the theory developed by Klein, although the point of view is different. Whereas in Klein's (and Cayley's) construction, the quadric is at infinity, the geometry, in the case developed by Poincar\'e, \emph{lives} on the quadric. 
 
 Let us recall that a \emph{quadric}, also called (by Poincar\'e) a \emph{quadratic surface},\index{quadratic surface} (``surface quadratique") is a surface in Euclidean three-space which is the zero locus of a  degree-two polynomial equation in three variables. There is a projective characterization of quadrics, which is coordinate-free: a quadric is a surface in projective space whose plane sections are all conics (real or imaginary). It follows from this definition that the intersection of any line with a quadric consists of two points, which may be real or imaginary, unless the line belongs to the quadric. Furthermore,  the set of all tangents to te quadrics from an arbitray point in space is a cone which cuts every plane in a conic, and the set of contact points of this cone with the quadric is also a conic.

 There are well-known classification of quadrics; some of them use coordinates and others are coordinate-free. Chapter 1 of the beautiful book of Hilbert and Cohn-Vossen, \emph{Geometry and the imagination} \cite{HCV}, concerns quadrics.  The equation of a conic can be put into normal form. Like for conics (which are the one-dimensional analogues of quadrics), there are some nondegenerate cases, and some degenerate cases. Poincar\'e obtained the two non-Euclidean geometries as geometries living on non-degenerate quadrics, and Euclidean geometry as a geometry living on a degenerate one. This is very close  to the ideas of Klein.

 There are nine types of quadrics. Six of them are \emph{ruled surfaces} (each point is on at least one straight line contained in the surface); these are  the cone, the one-sheeted hyperboloid, the hyperbolic paraboloid 
 and the three kind of cylinders (the elliptic, parabolic and hyperbolic). The three non-ruled quadrics are the ellipsoid, the elliptic paraboloid and the two-sheeted hyperboloid. These three surfaces do not contain any line.
 
The one-sheeted hyperboloid and the hyperbolic paraboloid are, like the plane, \emph{doubly ruled}, that is, each point is on at least two straight lines.

 Three types of nondegenerate quadrics which possess a center:  
 the ellipsoid, the two-sheeted hyperboloid and the one-sheeted hyperboloid.

Poincar\'e starts with a quadric in $\mathbb{R}^3$, called the \emph{fundamental surface}. On such a surface, he defines the notions of \emph{line}, of \emph{angle} between two lines and of \emph{length} of a segment.

 Given a quadric, the locus of midpoints of the system of chords that have a fixed direction is a plane, called a \emph{diametral plane} of the quadric. In the case where the quadric has a center, a diametral plane is a plane passing through the center. 
 
  Like in the case of the sphere, with its great circles and its small circles, Poincar\'e calls a \emph{line} an intersection of a quadric surface with a diametral plane, and a \emph{circle} an intersection of a quadric with an arbitrary plane. 
 
 Angles are then defined using the cross ratio. Given two lines $l_1$ and $l_2$ passing through a point $P$, Poincar\'e considers the quadruple of Euclidean lines formed by the tangents to $l_1$ and $l_2$ and the two rectilinear generatrices of the surface that pass through the point $P$. There are two generatrices at every point of the surface, and they may be real or imaginary.  Poincar\'e defines the \emph{angle} between $l_1$ and $l_2$ as the logarithm of the cross ratio of the four Euclidean lines ($l_1$, $l_2$ and the two generatrices) in the case where the two generatrices are real (and this occurs of the surface is a one-sheeted hyperboloid), and this logarithm divided by $\sqrt{-1}$ in the case where the generatrices are imaginary.
 
 Now given an arc of a line of the quadric, consider the cross ratio of the quadruple formed by the two extremities of this arc and the two points at infinity of the conic. The \emph{length} of this arc is the logarithm of the cross ratio of this quadruple of points in the case where the conic is a hyperbola, and the  logarithm of this cross ratio divided by $\sqrt{-1}$ otherwise.

 Poincar\'e then says that there are relations between  lengths and distances defined in this way, and that such relations constitute a set of theorems which are analogous to those of plane geometry. He calls the collection of theorems associated to a given quadric a \emph{quadratic geometry}. There are as many quadratic geometries as there are kinds of second degree surfaces, and Poincar\'e goes on with a classification of such geometries. 
 
 In the case where the fundamental surface is an ellipsoid, the geometry obtained is  spherical geometry.\index{geometry!spherical}
 
 In the case where the fundamental surface is a  two-sheeted hyperboloid, the geometry obtained is the Lobachevsky (or hyperbolic) geometry.\index{geometry!hyperbolic}
  
 In the case where the fundamental surface is an elliptic paraboloid, the geometry obtained is the Euclidean.\index{geometry!Euclidean} and Poincar\'e says that this geometry is a limiting geometry of each of  the previous two.
   
   There are other geometries, e.g. the one-sheeted hyperboloid and its various degenerate cases. Some of the degenerate geometries give the Euclidean geometry. But the one-sheeted hyperbolid itself gives a geometry which Poincar\'e highlights, as being a geometry which was not been studied yet, and  in which the following three phenomena occur:
   
1)   The distance between two points on the fundamental surface which are on a common rectilinear generatrix is zero.

2) There are two sorts of lines, those of the first kind, which correspond to the elliptic diametral sections, and those of the second kind, which correspond to the hyperbolic diametral sections. It is not possible, by a real motion, to bing a line of the first kind onto a line of the second kind.

3) There is no nontrivial real symmetry which sends a line onto itself. (Such a symmetry is possible in Euclidean geometry; it is obtained by a 180${}^{\mathrm{o}}$ rotation centered at a point on the line.)

This geometry is in fact the one called today the planar \emph{de Sitter geometry}.\index{geometry!de Sitter}

 Poincar\'e, in this paper,  does not mention Klein, but he thoroughly mentions Lie, and he considers this work as a consequence of Lie's work on groups. In the second part of the paper, titled ``Applications of group theory", Poincar\'e gives a characterization of the transformation group of each of these geometries. This is done in coordinates, at the infinitesimal level, in the tradition of Lie. He considers (p. 215) that ``geometry is nothing else than the study of a group".

 Our aim in this paper was to present to the reader of this book an important piece of work of Felix Klein. We also tried to  convey the idea that mathematical ideas occur at several people at the same period, when time is ready for that. Each  of us has a special way of thinking, and it often happens that works on the same problem, if they are not collective, complement themselves. We also hope that this paper will motivate the reader to go into the original sources.
 
 Let us conclude with the following two problems.
 
 \begin{enumerate}
 \item We already noted that Klein, in his development of the three geometries in his papers \cite{Klein-Ueber} and \cite{Klein-Ueber1}, considers Euclidean geometry as a \emph{transitional} geometry. In this way, Euclidean geometry corresponds to a limiting case of the absolute, in which the fundamental conic degenerates into a pair of imaginary points.
 We developed the notion of transitional geometry in the paper \cite{ACPK2}, in a way different from Klein's, and we studied  in which manner the fundamental notions of geometry (points, lines, distances, angles, etc.) as well as the trigonometric formulae transit from one geometry to another. An interesting problem is to make the same detailed study of transition of these fundamental notions in the context of Klein's description of the geometries.

 \item Hilbert developed a generalization of the Klein model of hyperbolic geometry where the underlying set is the interior of an ellipsoid to a geometry (called Hilbert geometry)\index{Hilbert geometry} where the underlying set is an arbitrary open convex set in $\mathbb{R}^n$.   The distance between two points $x$ and $y$ in Hilbert's generalization is again the logarithm of the cross ratio of the quadruple consisting of $x$ and $y$ and the two intersection points of the Euclidean line that joins these points with the boundary of the ellipsoid, taken in the natural order. We propose, as a problem, to develop generalizations of the two other geometries defined in the way Klein did it, that are analogous in some way to the generalization of hyperbolic geometry by Hilbert geometry.
  \end{enumerate}


\begin{thebibliography}{ABCD}
\bibitem{Taylor} B. Taylor, New principles of linear perspective, or the art of designing on a plane the representations of all sorts of objects, in a more general and simple method than has been done before, London, R. Knaplock, 1719.
\end{thebibliography}
\end{document}